\documentclass[leqno,draft]{article}

\newtheorem{theorem}{Theorem}
\newtheorem{lemma}[theorem]{Lemma}
\newtheorem{proposition}[theorem]{Proposition}
\newtheorem{definition}[theorem]{Definition}
\newtheorem{corollary}[theorem]{Corollary}

\newcommand{\begintheorem}{\addtocounter{equation}{1}\begin{theorem}}
\newcommand{\beginlemma}{\addtocounter{equation}{1}\begin{lemma}}
\newcommand{\beginproposition}{\addtocounter{equation}{1}\begin{proposition}}
\newcommand{\begindefinition}{\addtocounter{equation}{1}\begin{definition}}
\newcommand{\begincorollary}{\addtocounter{equation}{1}\begin{corollary}}

\begin{document}

\title{Groups and Cantor sets}

\author{Stephen Semmes \\
        Rice University}

\date{}

\maketitle

\begin{abstract}
It is well known that $n \times n$ upper-triangular real matrices with
$1$'s on the diagonal form a nilpotent Lie group with an interesting
family of non-isotropic dilations and corresponding geometry, as in
\cite{f-s}.  Here we look at $p$-adic versions of this, and related
matters.
\end{abstract}

\tableofcontents

\section{Semimetrics}
\label{semimetrics}
\setcounter{equation}{0}

        Let $X$ be a set.  A nonnegative real-valued function $d(x, y)$
on $X \times X$ is said to be a \emph{semimetric} on $X$, or equivalently
a \emph{pseudometric} on $X$, if $d(x, x) = 0$ for every $x \in X$,
\begin{equation}
\label{d(x, y) = d(y, x)}
        d(x, y) = d(y, x)
\end{equation}
for every $x, y \in X$, and
\begin{equation}
\label{d(x, z) le d(x, y) + d(y, z)}
        d(x, z) \le d(x, y) + d(y, z)
\end{equation}
for every $x, y, z \in X$.  If in addition $d(x, y) = 0$ only when $x
= y$, then $d(x, y)$ is a \emph{metric} on $X$.  The open ball in $X$
centered at a point $x \in X$ with radius $r > 0$ corresponding to a
semimetric $d(\cdot, \cdot)$ on $X$ is defined as usual by
\begin{equation}
\label{B(x, r) = {y in X : d(x, y) < r}}
        B(x, r) = \{y \in X : d(x, y) < r\}.
\end{equation}
Similarly, the closed ball in $X$ centered at $x \in X$ with radius $r
\ge 0$ corresponding to $d(\cdot, \cdot)$ is defined by
\begin{equation}
\label{overline{B}(x, r) = {y in X : d(x, y) le r}}
        \overline{B}(x, r) = \{y \in X : d(x, y) \le r\}.
\end{equation}

        Suppose that $V$ is a vector space over the real or complex numbers.
A nonnegative real-valued function $N(v)$ on $V$ is said to be a
\emph{seminorm} on $V$ if
\begin{equation}
\label{N(t v) = |t| N(v)}
        N(t \, v) = |t| \, N(v)
\end{equation}
for every $v \in V$ and $t \in {\bf R}$ or ${\bf C}$, as appropriate, and
\begin{equation}
\label{N(v + w) le N(v) + N(w)}
        N(v + w) \le N(v) + N(w)
\end{equation}
for every $v, w \in V$.  Here $|t|$ denotes the absolute value of $t$
in the real case, and the modulus of $t$ in the complex case.  If
$N(v) > 0$ when $v \ne 0$, then $N(v)$ is a \emph{norm} on $V$.  Thus
\begin{equation}
\label{d(v, w) = N(v - w)}
        d(v, w) = N(v - w)
\end{equation}
is a semimetric on $V$ when $V$ is a seminorm on $V$, and is a metric
on $V$ when $N(v)$ is a norm on $V$.

        If $d(x, y)$ is a metric on a set $X$, then there is a standard
Hausdorff topology on $X$ associated to $d(x, y)$, in which the open
balls in $X$ with respect to $d(x, y)$ are open sets, and the
collection of all such open balls forms a base for the topoology on
$X$.  In the same way, a semimetric on $X$ defines a topology on $X$
that may not be Hausdorff, which is generated by the corresponding
collection of open balls.  Suppose now that $I$ is a nonempty set, and
that $d_j(x, y)$ is a semimetric on $X$ for each $j \in I$.  Let
\begin{equation}
\label{B_j(x, r) = {y in X : d_j(x, y) < r}}
        B_j(x, r) = \{y \in X : d_j(x, y) < r\}
\end{equation}
be the open ball in $X$ centered at $x \in X$ with radius $r > 0$
corresponding to $d_j(\cdot, \cdot)$ for each $j \in I$.  One can
define a topology on $X$ by saying that a set $U \subseteq X$ is an
open set if for each $x \in U$ there are finitely many elements $j_1,
\ldots, j_n$ of $I$ and positive real numbers $r_1, \ldots, r_n$ such that
\begin{equation}
\label{bigcap_{k = 1}^n B_{j_k}(x, r_k) subseteq U}
        \bigcap_{k = 1}^n B_{j_k}(x, r_k) \subseteq U.
\end{equation}
It is easy to see that this is a topology on $X$, and that the open
balls associated to $d_j(x, y)$ are open sets in this topology for
each $j \in I$.  In order for this topology on $X$ to be Hausdorff, it
suffices that for each $x, y \in X$ with $x \ne y$ there be a $j \in
I$ such that $d_j(x, y) > 0$.

        If $d(x, y)$ is a semimetric on a set $X$ and $t$ is a positive
real number, then
\begin{equation}
\label{d'(x, y) = min(d(x, y), t)}
        d'(x, y) = \min(d(x, y), t)
\end{equation}
is also a semimetric on $X$.  Of course, this is a metric on $X$ when
$d(x, y)$ is a metric on $X$.  Let $B(x, r)$ be the open ball in $X$
centered at $x \in X$ with radius $r > 0$ associated to $d(\cdot,
\cdot)$, as in (\ref{B(x, r) = {y in X : d(x, y) < r}}), and let
\begin{equation}
\label{B'(x, r) = {y in X : d'(x, y) < r}}
        B'(x, r) = \{y \in X : d'(x, y) < r\}
\end{equation}
be the open ball corresponding to $d'(\cdot, \cdot)$.  By construction,
\begin{equation}
\label{B'(x, r) = B(x, r)}
        B'(x, r) = B(x, r)
\end{equation}
for every $x \in X$ and $r > 0$ such that $r \le t$, while
\begin{equation}
\label{B'(x, r) = X}
        B'(x, r) = X
\end{equation}
for every $x \in X$ when $r > t$.  It follows that the topology on $X$
determined by $d'(x, y)$ is the same as the topology determined by
$d(x, y)$.

        Similarly, let $I$ be a nonempty set, and let $d_j(x, y)$ be a 
semimetric on $X$ for each $j \in I$.  If $t_j$ is a positive real
number for each $j \in I$, then
\begin{equation}
\label{d_j'(x, y) = min(d_j(x, y), t_j)}
        d_j'(x, y) = \min(d_j(x, y), t_j)
\end{equation}
is a semimetric on $X$ for each $j \in I$, as before.  The topology on
$X$ determined by the semimetrics $d_j(x, y)$, $j \in I$, is the same
as the topology determined by the semimetrics $d_j'(x, y)$, $j \in I$.
One can also multiply semimetrics by positive real numbers, without
affecting the resulting topologies.

        Suppose that $I$ is a nonempty finite set, and let $d_j(x, y)$ be
a semimetric on a set $X$ for each $j \in I$.  Under these conditions,
\begin{equation}
\label{d(x, y) = max_{j in I} d_j(x, y)}
        d(x, y) = \max_{j \in I} d_j(x, y)
\end{equation}
is a semimetric on $X$, and the topology on $X$ associated to $d(x,
y)$ is the same as the topology determined by the collection of
semimetrics $d_j(x, y)$, $j \in I$.  More precisely, if $B(x, r)$ and
$B_j(x, r)$ are as in (\ref{B(x, r) = {y in X : d(x, y) < r}}) and
(\ref{B_j(x, r) = {y in X : d_j(x, y) < r}}), respectively, then
\begin{equation}
\label{B(x, r) = bigcap_{j in I} B_j(x, r)}
        B(x, r) = \bigcap_{j \in I} B_j(x, r)
\end{equation}
for each $x \in X$ and $r > 0$.  Note that $d(x, y)$ is a metric on
$X$ when for each $x, y \in X$ with $x \ne y$ there is a $j \in I$
such that $d_j(x, y) > 0$.

        Now let $d_1(x, y), d_2(x, y), d_3(x, y), \ldots$ be an infinite
sequence of semimetrics on a set $X$.  Also let $t_1, t_2, t_3,
\ldots$ be a sequence of positive real numbers that converges to $0$,
and let $d_j'(x, y)$ be the semimetric on $X$ given in (\ref{d_j'(x,
  y) = min(d_j(x, y), t_j)}) for each positive integer $j$.  Put
\begin{equation}
\label{d(x, y) = max_{j ge 1} d_j'(x, y)}
        d(x, y) = \max_{j \ge 1} d_j'(x, y)
\end{equation}
for each $x, y \in X$, which is equal to $0$ when $d_j'(x, y) = 0$ for
each $j$, or equivalently when $d_j(x, y) = 0$ for each $j$.
Otherwise, if $d_l'(x, y) > 0$ for some $l$, then
\begin{equation}
\label{d_j'(x, y) le t_j < d_l'(x, y)}
        d_j'(x, y) \le t_j < d_l'(x, y)
\end{equation}
for all sufficiently large $j$, which implies that the maximum in
(\ref{d(x, y) = max_{j ge 1} d_j'(x, y)}) is attained.  It is easy to
see that (\ref{d(x, y) = max_{j ge 1} d_j'(x, y)}) is a semimetric on
$X$, which is a metric on $X$ when for each $x, y \in X$ with $x \ne
y$ there is a $j \ge 1$ such that $d_j'(x, y) > 0$.

        Let $B(x, r)$ be the open ball in $X$ with center $x \in X$ and
radius $r > 0$ associated to (\ref{d(x, y) = max_{j ge 1} d_j'(x, y)}),
as in (\ref{B(x, r) = {y in X : d(x, y) < r}}), and let
\begin{equation}
\label{B_j'(x, r) = {y in X : d_j'(x, y) < r}}
        B_j'(x, r) = \{y \in X : d_j'(x, y) < r\}
\end{equation}
be the open ball corresponding to $d_j'(\cdot, \cdot)$ for each $j \ge
1$.  By construction,
\begin{equation}
\label{B(x, r) = bigcap_{j = 1}^infty B_j'(x, r)}
        B(x, r) = \bigcap_{j = 1}^\infty B_j'(x, r)
\end{equation}
for each $x \in X$ and $r > 0$, and $B_j'(x, r) = X$ when $t_j < r$.
It follows that
\begin{equation}
\label{B(x, r) = bigcap_{j = 1}^n B_j'(x, r)}
        B(x, r) = \bigcap_{j = 1}^n B_j'(x, r)
\end{equation}
when $n$ is sufficiently large so that $t_j < r$ for each $j > n$.
Using this, one can check that the topology on $X$ determined by
(\ref{d(x, y) = max_{j ge 1} d_j'(x, y)}) is the same as the topology
determined by the sequence of semimetrics $d_j'(x, y)$, which is the
same as the topology on $X$ determined by the original sequence of
semimetrics $d_j(x, y)$.

        Let $I$ be a nonempty set again, and let $d_j(\cdot, \cdot)$
be a semimetric on a set $X_j$ for each $j \in I$.  Also let $X = 
\prod_{j \in I} X_j$ be the Cartesian product of the $X_j$'s, $j \in I$.
An element of $X$ may be denoted $x = \{x_j\}_{j \in I}$, where $x_j \in X_j$
for each $j \in I$.  Put
\begin{equation}
\label{widetilde{d}_j(x, y) = d_j(x_j, y_j)}
        \widetilde{d}_j(x, y) = d_j(x_j, y_j)
\end{equation}
for each $x, y \in X$ and $j \in I$, which defines a semimetric on $X$
for each $j \in I$.  The topology on $X$ determined by the collection
of semimetrics $\widetilde{d}_j(x, y)$, $j \in I$, is the same as the
product topology on $X$ that corresponds to the topologies on the
$X_j$'s associated to the semimetrics $d_j(\cdot, \cdot)$.

        Let $a$, $b$ be nonnegative real numbers, and let $\alpha$ be
a positive real number less than or equal to $1$.  Observe that
\begin{equation}
\label{max(a, b) le (a^alpha + b^alpha)^{1/alpha}}
        \max(a, b) \le (a^\alpha + b^\alpha)^{1/\alpha},
\end{equation}
and hence that
\begin{eqnarray}
        a + b & \le & (a^\alpha + b^\alpha) \, (\max(a, b))^{1 - \alpha} \\
              & \le & (a^\alpha + b^\alpha)^{1 + (1 - \alpha)/\alpha}
                        = (a^\alpha + b^\alpha)^{1/\alpha}. \nonumber
\end{eqnarray}
This implies that
\begin{equation}
\label{(a + b)^alpha le a^alpha + b^alpha}
        (a + b)^\alpha \le a^\alpha + b^\alpha.
\end{equation}
If $d(x, y)$ is a semimetric on a set $X$, then it follows that $d(x,
y)^\alpha$ is also a semimetric on $X$ when $0 < \alpha \le 1$, which
is a metric on $X$ when $d(x, y)$ is a metric on $X$.  Note that $d(x,
y)^\alpha$ determines the same topology on $X$ as $d(x, y)$, and there
is an analogous statement for collections of semimetrics.

\section{Ultrametrics}
\label{ultrametrics}
\setcounter{equation}{0}

        A metric $d(x, y)$ on a set $X$ is said to be an \emph{ultrametric} if
\begin{equation}
\label{d(x, z) le max(d(x, y), d(y, z))}
        d(x, z) \le \max(d(x, y), d(y, z))
\end{equation}
for every $x, y, z \in X$, which automatically implies the usual
version (\ref{d(x, z) le d(x, y) + d(y, z)}) of the triangle
inequality.  Similarly, a semimetric $d(x, y)$ may be called a
\emph{semi-ultrametric} if it satisfies (\ref{d(x, z) le max(d(x, y),
  d(y, z))}) for every $x, y, z \in X$.  The discrete metric on any
set $X$ is an ultrametric, which is defined as usual by putting $d(x,
y) = 1$ when $x \ne y$ and $d(x, x) = 0$.  If $d(x, y)$ is any
semi-ultrametric on a set $X$ and $t > 0$, then (\ref{d'(x, y) =
  min(d(x, y), t)}) is a semi-ultrametric on $X$, as is $t \, d(x,
y)$.  It is easy to see that $d(x, y)^\alpha$ is also a
semi-ultrametric on $X$ for every $\alpha > 0$, which is a bit simpler
than the situation for ordinary semimetrics.

        Let $I$ be a nonempty set, and suppose that $d_j(x, y)$ is an
semi-ultrametric on a set $X$ for each $j \in I$.  If $I$ has only
finitely many elements, then the maximum (\ref{d(x, y) = max_{j in I}
  d_j(x, y)}) of $d_j(x, y)$ over $j \in I$ is a semi-ultrametric on
$X$ too.  Suppose instead that $I$ is the set ${\bf Z}_+$ of positive
integers, and that $t_1, t_2, t_3, \ldots$ is a sequence of positive
real numbers that converges to $0$.  In this case, (\ref{d(x, y) =
  max_{j ge 1} d_j'(x, y)}) is a semi-ultrametric on $X$ as well.  In
particular, this can be applied to the Cartesian product of a sequence
of sets, each of which is equipped with the discrete metric, and using
the corresponding sequence of semi-ultrametrics on the product, as in
(\ref{widetilde{d}_j(x, y) = d_j(x_j, y_j)}).

        Let $d(x, y)$ be a semimetric on a set $X$, and suppose that
$x, z \in X$, $r > 0$, and $d(x, y) < r$.  Thus $t = r - d(x, y) > 0$,
and it is easy to see that
\begin{equation}
\label{B(z, t) subseteq B(x, r)}
        B(z, t) \subseteq B(x, r),
\end{equation}
by the triangle inequality.  However, if $d(\cdot, \cdot)$ is a
semi-ultrametric on $X$, then (\ref{B(z, t) subseteq B(x, r)}) holds
with $t = r$.  Using this, one can check that $B(x, r)$ contains all
of its limit points in $X$, so that $B(x, r)$ is both open and closed
with respect to the topology on $X$ determined by $d(\cdot, \cdot)$.
Similarly,
\begin{equation}
\label{overline{B}(z, r) subseteq overline{B}(x, r)}
        \overline{B}(z, r) \subseteq \overline{B}(x, r)
\end{equation}
for every $x, z \in X$ and $r \ge 0$ such that $d(x, z) \le r$ when
$d(\cdot, \cdot)$ is a semi-ultrametric on $X$.  This implies that
$\overline{B}(x, r)$ is also both open and closed in $X$, with respect
to the topology determined by $d(\cdot, \cdot)$.  If $d(\cdot, \cdot)$
is an ultrametric on $X$, then it follows that $X$ is totally
disconnected with respect to the topology determined by $d(\cdot,
\cdot)$, in the sense that there are no connected subsets of $X$ with
more than one element.

        Let $X$ be a set, and let $\Delta$ be the diagonal subset of
$X \times X$, consisting of all ordered pairs $(x, x)$ with $x \in X$.
If $A$, $B$ are any subsets of $X \times X$, then put
\begin{equation}
\label{widetilde{A} = {(x, y) in X times X : (y, x) in A}}
        \widetilde{A} = \{(x, y) \in X \times X : (y, x) \in A\}
\end{equation}
and
\begin{eqnarray}
\label{A circ B = ...}
 A \circ B & = & \{(x, z) \in X \times X : \hbox{there is a } y \in X
                            \hbox{ such that } \\
 & & \qquad\qquad\qquad\qquad (x, y) \in A \hbox{ and } (y, z) \in B\}. 
                                                                 \nonumber
\end{eqnarray}
Suppose that $d(x, y)$ is a semimetric on $X$, and put
\begin{equation}
\label{U_r = {(x, y) in X times X : d(x, y) < r}}
        U_r = \{(x, y) \in X \times X : d(x, y) < r\}
\end{equation}
for each $r > 0$.  It is easy to see that
\begin{equation}
\label{Delta subseteq U_r = widetilde{U_r}}
        \Delta \subseteq U_r = \widetilde{U_r}
\end{equation}
for each $r > 0$,
\begin{equation}
\label{U_r subseteq U_t}
        U_r \subseteq U_t
\end{equation}
when $0 < r \le t$, and that
\begin{equation}
\label{U_r circ U_t subseteq U_{r + t}}
        U_r \circ U_t \subseteq U_{r + t}
\end{equation}
for every $r, t > 0$, by the triangle inequality.  If $d(x, y)$ is a
semi-ultrametric on $X$, then we have that
\begin{equation}
\label{U_r circ U_r subseteq U_r}
        U_r \circ U_r \subseteq U_r
\end{equation}
for every $r > 0$.

        As in \cite{jk}, a \emph{uniformity} on a set $X$ is a nonempty
collection $\mathcal{U}$ of subsets of $X \times X$ with the following
properties.  First, if $U, V \in \mathcal{U}$, then
\begin{equation}
\label{Delta subseteq U, widetilde{U} in mathcal{U}, and U cap V in mathcal{U}}
        \Delta \subseteq U, \ \widetilde{U} \in \mathcal{U}, \hbox{ and }
                                                 U \cap V \in \mathcal{U}.
\end{equation}
Second, for each $U \in \mathcal{U}$ there should be a $V \in \mathcal{U}$
such that
\begin{equation}
\label{V circ V subseteq U}
        V \circ V \subseteq U.
\end{equation}
Note that this implies that $V \subseteq U$, because $\Delta \subseteq
V$, as in (\ref{Delta subseteq U, widetilde{U} in mathcal{U}, and U
  cap V in mathcal{U}}).  Third, if $U \in \mathcal{U}$ and $U
\subseteq W \subseteq X \times X$, then $W \in \mathcal{U}$.  Thus $X
\times X \in \mathcal{U}$, because $\mathcal{U} \ne \emptyset$.  A set
$X$ with a uniformity $\mathcal{U}$ is said to be a \emph{uniform
  space}.  Roughly speaking, each element of $\mathcal{U}$ determines
a neighborhood around each point in $X$, which can be used to define
a topology on $X$ in particular.

        Let $d(x, y)$ be a semimetric on a set $X$, and let $U_r$ be as
in (\ref{U_r = {(x, y) in X times X : d(x, y) < r}}) for each $r > 0$.
It is easy to see that the collection $\mathcal{U}$ of subsets $U$ of 
$X \times X$ such that $U_r \subseteq U$ for some $r > 0$ is a uniformity
on $X$.  Similarly, let $I$ be a nonempty set, and suppose that
$d_j(x, y)$ is a semimetric on $X$ for each $j \in I$.  Put
\begin{equation}
\label{U_{r, j} = {(x, y) in X times X : d_j(x, y) < r}}
        U_{r, j} = \{(x, y) \in X \times X : d_j(x, y) < r\}
\end{equation}
for each $j \in I$ and $r > 0$, as in (\ref{U_r = {(x, y) in X times X
    : d(x, y) < r}}).  To get a uniformity on $X$, one can take the
collection of all subsets $U$ of $X \times X$ for which there are
finitely many elements $j_1, \ldots, j_n$ of $I$ and finitely many
positive real numbers $r_1, \ldots, r_n$ such that
\begin{equation}
\label{bigcap_{k = 1}^n U_{r_k, j_k} subseteq U}
        \bigcap_{k = 1}^n U_{r_k, j_k} \subseteq U.
\end{equation}
In both cases, the topology on $X$ associated to the uniformity is the
same as the one determined by the initial semimetric or collection of
semimetrices on $X$.  Conversely, it is well known that any uniformity
on $X$ corresponds to a collection of semimetrics on $X$ in this way.

        A subset $\mathcal{B}$ of a uniformity $\mathcal{U}$ on a set $X$
is said to be a \emph{base} for $\mathcal{U}$ if each element of
$\mathcal{U}$ contains an element of $\mathcal{B}$ as a subset.  If
$\mathcal{U}$ is the uniformity determined by a semimetric $d(x, y)$
on $X$ as before, then the collection of subsets $U_r$ of $X \times X$
with $r > 0$ is a base for $\mathcal{U}$.  Similarly, the collection
of subsets $U_r$ of $X \times X$ with $r = 1/n$ for some positive integer
$n$ is also a base for this uniformity.  Conversely, if there is a base
$\mathcal{B}$ for a uniformity $\mathcal{U}$ on $X$ such that $\mathcal{B}$
has only finitely or countably many elements, then it is well known that
there is a semimetric on $X$ for which $\mathcal{U}$ is the corresponding
uniformity.

        If $U$ is a subset of $X \times X$ for some set $X$ and $\Delta 
\subseteq U$, then
\begin{equation}
\label{U subseteq U circ U}
        U \subseteq U \circ U.
\end{equation}
In particular, if $U \subseteq X \times X$ satisfies
\begin{equation}
\label{Delta subseteq U and U circ U subseteq U}
        \Delta \subseteq U \hbox{ and } U \circ U \subseteq U,
\end{equation}
then $U \circ U = U$.  If $U, V \subseteq X \times X$ both satisfy
(\ref{Delta subseteq U and U circ U subseteq U}), then it is easy to
see that $\widetilde{U}$ and $U \cap V$ satisfy (\ref{Delta subseteq U
  and U circ U subseteq U}) too.  It follows that $W = U \cap
\widetilde{U}$ satisfies (\ref{Delta subseteq U and U circ U subseteq
  U}) when $U$ satisfies (\ref{Delta subseteq U and U circ U subseteq
  U}), and $\widetilde{W} = W$ automatically.  Note that $U \subseteq
X \times X$ satisfies
\begin{equation}
\label{Delta subseteq U, widetilde{U} = U, and U circ U subseteq U}
 \Delta \subseteq U, \ \widetilde{U} = U, \hbox{ and } U \circ U \subseteq U
\end{equation}
exactly when $U$ corresponds to an equivalence relation on $X$.  In
this case,
\begin{eqnarray}
\label{d_U(x, y) =  0 when (x, y) in U, ...}
 d_U(x, y) & = & 0 \quad\hbox{when } (x, y) \in U         \\
     & = & 1 \quad\hbox{when } (x, y) \in (X \times X) \setminus U. \nonumber
\end{eqnarray}
defines a semi-ultrametric on $X$.  Of course, if $d(x, y)$ is any
semi-ultrametric on $X$, then the set $U_r$ in (\ref{U_r = {(x, y) in
    X times X : d(x, y) < r}}) satisfies (\ref{Delta subseteq U,
  widetilde{U} = U, and U circ U subseteq U}) for each $r > 0$, by
(\ref{Delta subseteq U_r = widetilde{U_r}}) and (\ref{U_r circ U_r
  subseteq U_r}).

        Let $\mathcal{U}$ be a uniformity on $X$, and let $\mathcal{B}$
be a base for $\mathcal{U}$.  Thus $\Delta \subseteq U$ for every $U
\in \mathcal{B} \subseteq \mathcal{U}$, by definition of a uniformity.
If we also have that $\widetilde{U} = U$ and $U \circ U \subseteq U$
for each $U \in \mathcal{B}$, then $\mathcal{U}$ is the same as the
uniformity corresponding to the collection of semi-ultrametrics
$d_U(x, y)$ as in (\ref{d_U(x, y) = 0 when (x, y) in U, ...}) with $U
\in \mathcal{B}$.  It follows that $\mathcal{U}$ is the uniformity
corresponding to a single semi-ultrametric $d(x, y)$ when $\mathcal{U}$
satisfies this condition and has only finitely or countably many elements,
by the usual arguments.  If $\mathcal{B}$ is any base for $\mathcal{U}$, then
\begin{equation}
\label{mathcal{B}_1 = {U cap widetilde{U} : U in mathcal{B}}}
        \mathcal{B}_1 = \{U \cap \widetilde{U} : U \in \mathcal{B}\}
\end{equation}
is also a base for $\mathcal{B}$, and $\widetilde{W} = W$ for every $W
\in \mathcal{B}_1$.  If $\mathcal{B}$ is a base for $\mathcal{U}$ such
that $U \circ U \subseteq U$ for every $U \in \mathcal{B}$, then
$\widetilde{W} = W$ and $W \circ W \subseteq W$ for every $W \in
\mathcal{B}_1$, by the remarks in the preceding paragraph.  Hence the
previous observations about semi-ultrametrics can be applied to
$\mathcal{B}_1$ instead of $\mathcal{B}$.

\section{$p$-Adic numbers}
\label{p-adic numbers}
\setcounter{equation}{0}

        Let $p$ be a prime number, and let $|x|_p$ be the $p$-adic absolute 
value on the field ${\bf Q}$ of rational numbers.  Thus $|x|_p = 0$ when
$x = 0$, and otherwise $x$ can be expressed as $p^j \, a / b$ for some
integers $a$, $b$, and $j$, where $a$ and $b$ are nonzero and not integer
multiples of $p$.  In this case,
\begin{equation}
\label{|x|_p = p^{-j}}
        |x|_p = p^{-j}.
\end{equation}
It is well known and easy to see that
\begin{equation}
\label{|x + y|_p le max(|x|_p, |y|_p)}
        |x + y|_p \le \max(|x|_p, |y|_p)
\end{equation}
and
\begin{equation}
\label{|x y|_p = |x|_p |y|_p}
        |x \, y|_p = |x|_p \, |y|_p
\end{equation}
for every $x, y \in {\bf Q}$.  The $p$-adic metric is defined on ${\bf Q}$ by
\begin{equation}
\label{d_p(x, y) = |x - y|_p}
        d_p(x, y) = |x - y|_p,
\end{equation}
and this is an ultrametric on ${\bf Q}$, because of (\ref{|x + y|_p le
  max(|x|_p, |y|_p)}).

        The field ${\bf Q}_p$ of $p$-adic numbers can be obtained
by completing ${\bf Q}$ as a metric space with respect to the $p$-adic
metric, in essentially the same way that the field ${\bf R}$ of real
numbers can be obtained by completing ${\bf Q}$ with respect to the
standard metric.  In particular, the $p$-adic absolute value and metric
have natural extensions to ${\bf Q}_p$ with the same properties as before,
and the field operations on ${\bf Q}_p$ are continuous with respect to
the $p$-adic metric.  It follows that
\begin{equation}
\label{{bf Z}_p = {x in {bf Q}_p : |x|_p le 1}}
        {\bf Z}_p = \{x \in {\bf Q}_p : |x|_p \le 1\}
\end{equation}
is a closed sub-ring of ${\bf Q}_p$, known as the ring of $p$-adic
integers.  This contains the usual ring ${\bf Z}$ of integers by
definition of the $p$-adic absolute value, and one can show that ${\bf
  Z}_p$ is the same as the closure of ${\bf Z}$ in ${\bf Q}_p$.  Note
that $|x|_p$ is an integer power of $p$ for every $x \in {\bf Q}_p$
with $x \ne 0$.

        If $j \in {\bf Z}$, then put
\begin{equation}
\label{p^j Z_p = {p^j x : x in Z_p} = {y in Q_p : |y|_p le p^{-j}}}
        p^j \, {\bf Z}_p = \{p^j \, x : x \in {\bf Z}_p\} 
                         = \{y \in {\bf Q}_p : |y|_p \le p^{-j}\}.
\end{equation}
This is a closed subgroup of ${\bf Q}_p$ as a commutative group with
respect to addition for each $j \in {\bf Z}$, and an ideal in ${\bf
  Z}_p$ when $j \ge 0$.  If $j \ge 0$, then the inclusions ${\bf Z}
\subseteq {\bf Z}_p$ and $p^j \, {\bf Z} \subseteq p^j \, {\bf Z}_p$
lead to a ring homomorphism from ${\bf Z} / p^j \, {\bf Z}$ into ${\bf
  Z}_p / p^j \, {\bf Z}_p$, and one can check that this homomorphism
is an isomorphism.  In particular, ${\bf Z}_p / p^j \, {\bf Z}_p$ has
exactly $p^j$ elements for each nonnegative integer $j$, so that ${\bf
  Z}_p$ is the union of $p^j$ pairwise-disjoint translates of $p^j \,
{\bf Z}_p$.  This implies that ${\bf Z}_p$ is totally bounded in ${\bf
  Q}_p$, and hence that ${\bf Z}_p$ is compact, since ${\bf Z}_p$ is
also a closed set in ${\bf Q}_p$ and ${\bf Q}_p$ is compact.  It
follows that $p^l \, {\bf Z}_p$ is a compact set in ${\bf Q}_p$ for
every $l \in {\bf Z}$, and that closed and bounded sets in ${\bf Q}_p$
are compact.  If $H$ is Haar measure on ${\bf Q}_p$ normalized so that
$H({\bf Z}_p) = 1$, then
\begin{equation}
\label{H(p^l {bf Z}_p) = p^{-l}}
        H(p^l \, {\bf Z}_p) = p^{-l}
\end{equation}
for each $l \in {\bf Z}$.

\section{Ultranorms}
\label{ultranorms}
\setcounter{equation}{0}

        Let $p$ be a prime number, and let $V$ be a vector space
over the $p$-adic numbers ${\bf Q}_p$.  A nonnegative real-valued 
function $N(v)$ on $V$ is said to be an \emph{ultranorm} on $V$
if $N(v) > 0$ when $v \ne 0$,
\begin{equation}
\label{N(t v) = |t|_p N(v)}
        N(t \, v) = |t|_p \, N(v)
\end{equation}
for every $v \in V$ and $t \in {\bf Q}_p$, and
\begin{equation}
\label{N(v + w) le max(N(v), N(w))}
        N(v + w) \le \max(N(v), N(w))
\end{equation}
for every $v, w \in V$.  In this case, it is easy to see that
\begin{equation}
\label{d(v, w) = N(v - w), 2}
        d(v, w) = N(v - w)
\end{equation}
is an ultrametric on $V$, and that addition and scalar multiplication
on $V$ are continuous with respect to the topology determined by
(\ref{d(v, w) = N(v - w), 2}).

        Let $n$ be a positive integer, and let ${\bf Q}_p^n$ be the space
of $n$-tuples of $p$-adic numbers, which is a vector space over ${\bf Q}_p$
with respect to coordinatewise addition and scalar multiplication.  Put
\begin{equation}
\label{N(v) = max(|v_1|_p, ldots, |v_n|_p)}
        N(v) = \max(|v_1|_p, \ldots, |v_n|_p)
\end{equation}
for each $v = (v_1, \ldots, v_n) \in {\bf Q}_p^n$, which clearly
defines an ultranorm on ${\bf Q}_p^n$.  Note that the topology on
${\bf Q}_p^n$ determined by $N$ is the same as the product topology
corresponding to the topology on ${\bf Q}_p$ determined by the
$p$-adic metric.

        Similarly, let $M_n({\bf Q}_p)$ be the space of $n \times n$
matrices with entries in ${\bf Q}_p$, which is a vector space over
${\bf Q}_p$ with respect to entry-wise addition and scalar multiplication.
Put
\begin{equation}
\label{||A|| = max_{1 le j, k le n} |a_{j, k}|_p}
        \|A\| = \max_{1 \le j, k \le n} |a_{j, k}|_p
\end{equation}
for each $A = \{a_{j,k}\}_{j, k = 1}^n \in M_n({\bf Q}_p)$, which
defines an ultranorm on $M_n({\bf Q}_p)$.  Of course, we can identify
$M_n({\bf Q}_p)$ with ${\bf Q}_p^{n^2}$ in an obvious way, so that
(\ref{||A|| = max_{1 le j, k le n} |a_{j, k}|_p}) corresponds to the
standard ultranorm defined in the previous paragraph.  If $A = \{a_{j,
  k}\}_{j, k = 1}^n$ and $B = \{b_{j, k}\}_{j, k = 1}^n$ are elements
of $M_n({\bf Q}_p)$, then their product $A \, B = C = \{c_{j, l}\}_{j,
  l = 1}^n$ can be defined as usual by
\begin{equation}
\label{c_{j, l} = sum_{k = 1}^n a_{j, k} b_{k, l}}
        c_{j, l} = \sum_{k = 1}^n a_{j, k} \, b_{k, l}.
\end{equation}
Observe that
\begin{equation}
\label{|c_{j, l}|_p le max_{1 le k le n} (|a_{j, k}|_p |b_{k, l}|_p)}
        |c_{j, l}|_p \le \max_{1 \le k \le n} (|a_{j, k}|_p \, |b_{k, l}|_p)
\end{equation}
for each $j, l = 1, \ldots, n$, which implies that
\begin{equation}
\label{||A B|| le ||A|| ||B||}
        \|A \, B\| \le \|A\| \, \|B\|.
\end{equation}

        Let $\widehat{A}$ be the linear mapping from ${\bf Q}_p^n$ into
itself corresponding to $A = \{a_{j, k}\}_{j, k = 1}^n \in M_n({\bf Q}_p)$,
so that the $j$th component of $\widehat{A}(v)$ is equal to
\begin{equation}
\label{(widehat{A}(v))_j = sum_{k = 1}^n a_{j, k} v_k}
        (\widehat{A}(v))_j = \sum_{k = 1}^n a_{j, k} \, v_k
\end{equation}
for each $v \in {\bf Q}_p^n$ and $j = 1, \ldots, n$.  If $B = \{b_{j,
  k}\}_{j, k = 1}^n$ is another element of $M_n({\bf Q}_p)$, then
\begin{equation}
\label{widehat{(A B)}(v) = ... = widehat{A}(widehat{B}(v))}
        \widehat{(A \, B)}(v) = (\widehat{A} \circ \widehat{B})(v) 
                              = \widehat{A}(\widehat{B}(v))
\end{equation}
for every $v \in {\bf Q}_p^n$, so that matrix multiplication corresponds
to composition of linear mappings in the usual way.  It follows from
(\ref{(widehat{A}(v))_j = sum_{k = 1}^n a_{j, k} v_k}) that
\begin{equation}
\label{|(widehat{A}(v))_j|_p le max_{1 le k le n} |a_{j, k}|_p |v_k|_p}
        |(\widehat{A}(v))_j|_p \le \max_{1 \le k \le n} |a_{j, k}|_p \, |v_k|_p
\end{equation}
for each $v \in {\bf Q}_p^n$ and $j = 1, \ldots, n$, which implies that
\begin{equation}
\label{N(widehat{A}(v)) le ||A|| N(v)}
        N(\widehat{A}(v)) \le \|A\| \, N(v)
\end{equation}
for every $v \in {\bf Q}_p^n$, where $N(v)$ is as in (\ref{N(v) =
  max(|v_1|_p, ldots, |v_n|_p)}).  It is easy to see that $\|A\|$ is
the smallest nonnegative real number with this property, by
considering the case where $v$ is a standard basis vector for ${\bf
  Q}_p^n$, with one component equal to $1$ and the rest equal to $0$.
Thus $\|A\|$ is the same as the operator norm of $\widehat{A}$ with
respect to the standard ultranorm $N(v)$ on ${\bf Q}_p^n$.

\section{Topological groups}
\label{topological groups}
\setcounter{equation}{0}

        Let $G$ be a topological group, so that $G$ is both a group
and a topological space, and the group operations are continuous.
It is customary to ask also that the set consisting of the identity
element $e$ in $G$ be a closed set, which implies that $G$ is
Hausdorff as a topological space.  Using multiplicative notation
for the group operations, a semimetric $d(x, y)$ on $G$ is said to 
be invariant under left translations on $G$ if
\begin{equation}
\label{d(a x, a y) = d(x, y)}
        d(a \, x, a \, y) = d(x, y)
\end{equation}
for every $a, x, y \in G$.  Similarly, $d(x, y)$ is invariant under
right translations on $G$ if
\begin{equation}
\label{d(x b, y b) = d(x, y)}
        d(x \, b, y \, b) = d(x, y)
\end{equation}
for every $b, x, y \in G$.  Note that $d(x, y)$ is invariant under
left translations if and only if $d(x^{-1}, y^{-1})$ is invariant
under right translations, and that the discrete metric on $G$ is
invariant under both left and right translations.

        It is well known that there is a collection of left-invariant 
semimetrics on $G$ that determines the same topology for $G$.  If
there is a local base for the topology of $G$ at $e$ with only
finitely or countably many elements, then there is a left-invariant
metric on $G$ that determines the same topology.  More precisely, if
there is a local base for the topology of $G$ at $e$ with only
finitely many elements, then $\{e\}$ is an open set in $G$, so that
$G$ is equipped with the discrete topology, and one can simply use the
discrete metric on $G$.  One can define a left-invariant uniformity on
$G$ more directly which is compatible with the given topology, and
with any collection of left-invariant semimetrics on $G$ that
determines the same topology.  There are analogous statements for
right-invariant semimetrics, which can be derived from the previous
statements for left-invariant semimetrics using the mapping $x \mapsto
x^{-1}$.

        Let $d(x, y)$ be a left-invariant semimetric on a group $G$,
and put
\begin{equation}
\label{r(x) = d(x, e)}
        r(x) = d(x, e)
\end{equation}
for each $x \in G$.  Thus $r(x)$ is a nonnegative real-valued function
on $G$, $r(e) = 0$, and we can use the symmetry and left-invariance of
$d(x, y)$ to get that
\begin{equation}
\label{r(x) = d(x, e) = d(e, x) = d(x^{-1}, e) = r(x^{-1})}
        r(x) = d(x, e) = d(e, x) = d(x^{-1}, e) = r(x^{-1})
\end{equation}
for each $x \in G$.  Similarly,
\begin{equation}
\label{r(x y) = d(x y, e) le d(x y, x) + d(x, e) = ... = r(x) + r(y)}
 \quad  r(x \, y) = d(x \, y, e) \le d(x \, y, x) + d(x, e) 
                  = d(y, e) + d(x, e) = r(x) + r(y)
\end{equation}
for every $x, y \in G$, using the triangle inequality and
left-invariance.  If $d(x, y)$ is a left-invariant semi-ultrametric on
$G$, then
\begin{eqnarray}
\label{r(x y) = d(x y, e) le max(d(x y, x), d(x, e)) = ... = max(r(x), r(y))}
        r(x \, y) = d(x \, y, e) & \le & \max(d(x \, y, x), d(x, e)) \\
                     & = & \max(d(y, e), d(x, e)) = \max(r(x), r(y)) \nonumber
\end{eqnarray}
for every $x, y \in G$.

        Conversely, suppose that $r(x)$ is a nonnegative real-valued
function on $G$ such that $r(e) = 0$,
\begin{equation}
\label{r(x) = r(x^{-1})}
        r(x) = r(x^{-1})
\end{equation}
for every $x \in G$, and
\begin{equation}
\label{r(x y) le r(x) + r(y)}
        r(x \, y) \le r(x) + r(y)
\end{equation}
for every $x, y \in G$.  Under these conditions, one can check that
\begin{equation}
\label{d(x, y) = r(y^{-1} x)}
        d(x, y) = r(y^{-1} \, x)
\end{equation}
defines a left-invariant semimetric on $G$.  If $r(x)$ also satisfies
\begin{equation}
\label{r(x y) le max(r(x), r(y))}
        r(x \, y) \le \max(r(x), r(y))
\end{equation}
for every $x, y \in G$, then (\ref{d(x, y) = r(y^{-1} x)}) defines a
semi-ultrametric on $G$.  Of course, if $r(x)$ is as in (\ref{r(x) =
  d(x, e)}) for some left-invariant semimetric $d(x, y)$ on $G$,
then (\ref{d(x, y) = r(y^{-1} x)}) holds for every $x, y \in G$.

        If $d(x, y)$ is instead a right-invariant semimetric on $G$, then 
(\ref{r(x) = d(x, e)}) still satisfies $r(e) = 0$, (\ref{r(x) = r(x^{-1})}), 
and (\ref{r(x y) le r(x) + r(y)}).  This can be verified in the same
way as before, or using the fact that $d(x^{-1}, y^{-1})$ is a
left-invariant semimetric.  If $d(x, y)$ is a right-invariant
semi-ultrametric on $G$, then (\ref{r(x y) le max(r(x), r(y))}) holds
as well.  Conversely, if $r(x)$ is a nonnegative real-valued function
on $G$ that satisfies $r(e) = 0$, (\ref{r(x) = r(x^{-1})}), and
(\ref{r(x y) le r(x) + r(y)}), then
\begin{equation}
\label{d(x, y) = r(x y^{-1})}
        d(x, y) = r(x \, y^{-1})
\end{equation}
defines a right-invariant semimetric on $G$.  If $r(x)$ also satisfies
(\ref{r(x y) le max(r(x), r(y))}), then (\ref{d(x, y) = r(x y^{-1})})
is a semi-ultrametric on $G$.  As before, if $r(x)$ is as in
(\ref{r(x) = d(x, e)}) for some right-invariant semimetric $d(x, y)$
on $G$, then (\ref{d(x, y) = r(x y^{-1})}) holds for every $x, y \in
G$.  Note that (\ref{d(x, y) = r(y^{-1} x)}) and (\ref{d(x, y) = r(x
  y^{-1})}) correspond to each other under the mapping $x \mapsto
x^{-1}$, because of (\ref{r(x) = r(x^{-1})}).

        If $d(x, y)$ is a semimetric on $G$ that is invariant under
both left and right translations, then (\ref{r(x) = d(x, e)}) has the 
additional property that
\begin{equation}
\label{r(a x a^{-1}) = r(x)}
        r(a \, x \, a^{-1}) = r(x)
\end{equation}
for every $a, x \in G$.  This is equivalent to the condition that
\begin{equation}
\label{r(x y) = r(y x)}
        r(x \, y) = r(y \, x)
\end{equation}
for every $x, y \in G$.  Conversely, if $r(x)$ is a nonnegative
real-valued function on $G$ that satisfies $r(e) = 0$, (\ref{r(x) =
  r(x^{-1})}), (\ref{r(x y) le r(x) + r(y)}), and (\ref{r(x y) = r(y
  x)}), then (\ref{d(x, y) = r(y^{-1} x)}) and (\ref{d(x, y) = r(x
  y^{-1})}) are equal to each other, and hence define a semimetric on
$G$ that is invariant under both left and right translations.

        Let us say that a semimetric $d(x, y)$ on a topological group
$G$ is compatible with the topology on $G$ if every open ball in $G$
with respect to $d(x, y)$ is an open set in $G$.  This is equivalent
to asking that each open ball in $G$ with respect to $d(x, y)$ contain
an open set in $G$ that contains the center of the ball.  If $d(x, y)$
is invariant under left or right translations on $G$, then it suffices
to check this condition for balls centered at $e$.  Alternatively,
this means that every open set in $G$ with respect to the topology
determined by $d(x, y)$ is also an open set with respect to the given
topology on $G$.  In particular, if the topology on $G$ is determined
by a collection of semimetrics, then each of the semimetrics in the
collection is compatible with the topology on $G$ in this sense.

\section{Small open subgroups}
\label{small open subgroups}
\setcounter{equation}{0}

        Let $G$ be a group, and let $d(x, y)$ be a semi-ultrametric on $G$
that is invariant under left or right translations.  Also let $r(x)$
be as in (\ref{r(x) = d(x, e)}), so that $r(x)$ is a nonnegative real-valued
function on $G$ that satisfies $r(e) = 0$, (\ref{r(x) = r(x^{-1})}),  and 
(\ref{r(x y) le max(r(x), r(y))}).  If $t$ is a positive real number, then
it follows that
\begin{equation}
\label{U_t = {x in G : r(x) < t}}
        U_t = \{x \in G : r(x) < t\}
\end{equation}
is a subgroup of $G$.  This is the same as the open ball in $G$
centered at $e$ with radius $t$ with respect to $d(x, y)$, which is an
open set in $G$ when $G$ is a topological group and $d(x, y)$ is
compatible with the topology on $G$.  If $d(x, y)$ is invariant under
both left and right translations on $G$, then $r(x)$ satisfies
(\ref{r(a x a^{-1}) = r(x)}), and hence (\ref{U_t = {x in G : r(x) <
    t}}) is a normal subgroup in $G$ for each $t > 0$.

        Suppose that $U$ is a subgroup of $G$, and let $r_U(x)$ be the 
indicator function associated to $U$ on $G$, which is equal to $0$
when $x \in U$ and to $1$ when $x \in G$ is not in $U$.  Thus $r_U(x)$
satisfies $r_U(e) = 0$, (\ref{r(x) = r(x^{-1})}), and (\ref{r(x y) le
  max(r(x), r(y))}), which implies that (\ref{d(x, y) = r(y^{-1} x)}) and 
(\ref{d(x, y) = r(x y^{-1})}) are semi-ultrametrics on $G$ that are invariant
under left and right translations, respectively.  If $G$ is a topological
group and $U$ is an open subgroup of $G$, then these semi-ultrametrics
are compatible with the topology on $U$.  If $U$ is a normal subgroup
of $G$, then (\ref{d(x, y) = r(y^{-1} x)}) and (\ref{d(x, y) = r(x y^{-1})})
are equal to each other, and hence define a semi-ultrametric on $G$
that is invariant under both left and right translations.

        Let us say that a topological group $G$ has \emph{small open 
subgroups} if there is a local base for the topology of $G$ at $e$ 
consisting of open subgroups of $G$.  This happens if and only if
there is a collection of semi-ultrametrics on $G$ that are invariant
under left or right translations and which determine the same topology
on $G$, by the remarks in the previous paragraphs.  If there is a
countable local base for the topology of $G$ consisting of open
subgroups, then there is an ultrametric on $G$ that determines the
same topology on $G$ and which is invariant under left translations,
or under right translations.  If $G$ has small open subgroups, then
$G$ is totally disconnected, and in fact has topological dimension
$0$.  Conversely, if $G$ is locally compact and totally disconnected,
then it is well known that $G$ has small open subgroups.

        Similarly, let us say that $G$ has \emph{small open normal 
subgroups} if there is a local base for the topology of $G$ consisting 
of open normal subgroups of $G$.  As before, this happens if and only
if there is a collection of semi-ultrametrics on $G$ that are
invariant under both left and right translations and which determine
the same topology on $G$, by the remarks at the beginning of the
section.  If there is a countable local base for the topology of $G$
consisting of open normal subgroups of $G$, then there is an
ultrametric on $G$ that is invariant under both left and right
translations and which determines the same topology on $G$.

        If $G$ is any group equipped with the discrete topology, then 
$\{e\}$ is an open normal subgroup of $G$, so that $G$ has small open
subgroups in particular.  If $G$ is a product of discrete groups with
the product topology, then $G$ has small open normal subgroups, and
the discrete metrics on the factors lead to semimetrics on $G$, as in
(\ref{widetilde{d}_j(x, y) = d_j(x_j, y_j)}).  Each of these
semimetrics is invariant under both left and right translations on
$G$, and the collection of these semimetrics determines the product
topology on $G$, as before.  Of course, the product of finitely many
discrete groups is also discrete.  If $G$ is the product of countably
many nontrivial discrete groups, then there is a countable local base
for the topology of $G$ at $e$ consisting open normal subgroups of $G$, 
and there is an ultrametric on $G$ that is invariant under both left and
right translations and which determines the product topology on $G$.

        Let $p$ be a prime number, and consider ${\bf Q}_p$ as a commutative 
topological group with respect to addition.  Clearly ${\bf Q}_p$ has
small open subgroups, because $p^j \, {\bf Z}_p$ is an open subgroup
of ${\bf Q}_p$ for each integer $j$.  The $p$-adic metric on ${\bf
  Q}_p$ is a traslation-invariant ultrametric on ${\bf Q}_p$.

        The group ${\bf Q}_p^*$ of nonzero $p$-adic numbers with respect to 
multiplication is also a commutative topological group, and
\begin{equation}
\label{{x in {bf Q}_p : |x|_p = 1}}
        \{x \in {\bf Q}_p : |x|_p = 1\}
\end{equation}
is an open subgroup of ${\bf Q}_p^*$.  Similarly,
\begin{equation}
\label{{x in {bf Q}_p : x - 1 in p^j {bf Z}_p}}
        \{x \in {\bf Q}_p : x - 1 \in p^j \, {\bf Z}_p\}
\end{equation}
is an open subgroup of ${\bf Q}_p^*$ for each positive integer $j$, so
that ${\bf Q}_p^*$ has small open subgroups too.

        Observe that the $p$-adic metric (\ref{d_p(x, y) = |x - y|_p}) 
satisfies
\begin{equation}
\label{d_p(a x, a y) = |a x - a y|_p = |x - y|_p}
        d_p(a \, x, a \, y) = |a \, x - a \, y|_p = |x - y|_p
\end{equation}
for every $a, x, y \in {\bf Q}_p$ with $|a|_p = 1$.  In particular,
the restriction of $d_p(x, y)$ to $x, y \in {\bf Q}_p$ with $|x|_p =
|y|_p = 1$ is a translation-invariant ultrametric on (\ref{{x in {bf
      Q}_p : |x|_p = 1}}) as a compact commutative topological group
with respect to multiplication.

        Put
\begin{equation}
\label{r_p(x) = max(|x - 1|_p, |(1/x) - 1|_p)}
        r_p(x) = \max(|x - 1|_p, |(1/x) - 1|_p)
\end{equation}
for each $x \in {\bf Q}_p^*$, so that $r_p(1/x) = r_p(x)$ by
construction.  Equivalently,
\begin{equation}
\label{r_p(x) = |x - 1|_p = |(1/x) - 1|_p le 1}
        r_p(x) = |x - 1|_p = |(1/x) - 1|_p \le 1
\end{equation}
when $|x|_p = 1$, and otherwise
\begin{equation}
\label{r_p(x) = max(|x|_p, 1/|x|_p) ge p}
        r_p(x) = \max(|x|_p, 1/|x|_p) \ge p
\end{equation}
when $|x|_p \ne 1$.  If $x, y \in {\bf Q}_p^*$ satisfy $|x|_p = |y|_p
= 1$, then
\begin{eqnarray}
\label{r_p(x y) = |x y - 1|_p le ... = max(r_p(x), r_p(y))}
 r_p(x \, y) = |x \, y - 1|_p & \le & \max(|x \, y - y|_p, |y - 1|_p) \\
            & = & \max(|x - 1|_p, |y - 1|_p) = \max(r_p(x), r_p(y)). \nonumber
\end{eqnarray}
Note that $r_p(1) = 0$, and that $r_p(x) > 0$ when $x \ne 1$.

        Let $t$ be a real number such that $1 \le t \le p$, and put
\begin{equation}
\label{r_p'(x) = min(r_p(x), t)}
        r_p'(x) = \min(r_p(x), t)
\end{equation}
for each $x \in {\bf Q}_p^*$.  Thus
\begin{equation}
\label{r_p'(x) = r_p(x)}
        r_p'(x) = r_p(x)
\end{equation}
when $|x|_p = 1$, and
\begin{equation}
\label{r_p'(x) = t}
        r_p'(x) = t
\end{equation}
when $|x|_p \ne 1$.  Let us check that
\begin{equation}
\label{r_p'(x y) le max(r_p'(x), r_p'(y))}
        r_p'(x \, y) \le \max(r_p'(x), r_p'(y))
\end{equation}
for every $x, y \in {\bf Q}_p^*$.  If $|x|_p = |y|_p = 1$, then this
reduces to (\ref{r_p(x y) = |x y - 1|_p le ... = max(r_p(x),
  r_p(y))}), by (\ref{r_p'(x) = r_p(x)}).  Otherwise, if $|x|_p \ne 1$
or $|y|_p \ne 1$, then the right side of (\ref{r_p'(x y) le
  max(r_p'(x), r_p'(y))}) is equal to $t$, so that (\ref{r_p'(x y) le
  max(r_p'(x), r_p'(y))}) holds trivially.  Of course, $r_p'(1/x) =
r_p'(x)$ for each $x \in {\bf Q}_p^*$, because of the analogous
property for $r_p(x)$.  Similarly, $r_p'(1) = 0$, and $r_p'(x) > 0$
for each $x \in {\bf Q}_p^*$ with $x \ne 1$.

        It follows that
\begin{equation}
\label{d_p'(x, y) = r_p'(x / y)}
        d_p'(x, y) = r_p'(x / y)
\end{equation}
defines a translation-invariant ultrametric on ${\bf Q}_p^*$, as a
commutative group with respect to multiplication.  If $|x|_p = |y|_p$,
then $|x / y|_p = 1$, and hence
\begin{equation}
\label{d_p'(x, y) = r_p'(x / y) = r_p(x / y) = ... = |x - y|_p / |y|_p le 1}
        d_p'(x, y) = r_p'(x / y) = r_p(x / y) = |(x / y) - 1|_p 
                                              = \frac{|x - y|_p}{|y|_p} \le 1.
\end{equation}
Otheriwse, if $|x|_p \ne |y|_p$, then $d_p'(x, y) = t$, by
(\ref{r_p'(x) = t}).  Note that $|x|_p = |y|_p$ when $|x - y|_p <
\max(|x|_p, |y|_p)$, because of the ultrametric version of the
triangle inequality.  This shows that the topology on ${\bf Q}_p^*$
determined by $d_p'(x, y)$ is the same as the topology induced by the
usual one on ${\bf Q}_p$.

\section{Large compact subgroups}
\label{large compact subgroups}
\setcounter{equation}{0}

        Let us say that a topological group $G$ has \emph{large compact
subgroups} if every compact set in $G$ is contained in a compact subgroup
of $G$.  Similarly, $G$ has \emph{large compact open subgroups} if every
compact set in $G$ is contained in a compact open subgroup of $G$.
If a subgroup $H$ of $G$ contains a nonempty open set, then $H$ is
an open set as well.  In particular, if a locally compact topological
group $G$ has large compact subgroups, then it automatically has large
compact open subgroups.  Of course, there are analogous conditions using
compact normal subgroups, and compact topological groups have all 
of these properties trivially.

        Observe that ${\bf Q}_p$ has large compact subgroups as a 
commutative topological group with respect to addition for each prime
number $p$.  More precisely, $p^j \, {\bf Z}_p$ is a compact open subgroup
of ${\bf Q}_p$ for each integer $j$, and every compact set in ${\bf Q}_p$
is contained in $p^j \, {\bf Z}_p$ for some $j$.  However, the multiplicative
group ${\bf Q}_p^*$ of nonzero complex numbers does not have large compact
subgroups.  Indeed, if $x \in {\bf Q}_p^*$ and $|x|_p \ne 1$, then there is no
compact subgroup of ${\bf Q}_p^*$ that contains $x$.  The additive group
${\bf Z}$ of integers does not have large compact subgroups with respect
to the discrete topology, since $\{0\}$ is the only compact subgroup of
${\bf Z}$.

        Let $G$ be a topological group, and let $d(x, y)$ be a semimetric on
$G$ which is compatible with the topology on $G$.  Note that
\begin{equation}
\label{|d(x, y) - d(x, z)| le d(y, z)}
        |d(x, y) - d(x, z)| \le d(y, z)
\end{equation}
for every $x, y, z \in G$, because of the triangle inequality.  This
implies that for each $x \in G$, $d(x, y)$ is a continuous function of
$y$ on $G$, when $d(x, y)$ is compatible with the topology on $G$.  It
follows that for each $y \in G$, $d(x, y)$ is a continuous function of
$x$ on $G$, and one can also check that $d(x, y)$ is continuous as a
function of $x$ and $y$ on $G \times G$.  In particular, compact
subsets of $G$ are bounded with respect to $d(x, y)$ under these
conditions.

        Let us say that $d(x, y)$ is \emph{proper} on $G$ if closed
subsets of $G$ that are bounded with respect to $d(x, y)$ are compact
with respect to the topology on $G$.  Of course, closed balls with
respect to $d(x, y)$ are closed sets in $G$ when $d(x, y)$ is compatible
with the topology on $G$, and hence they are compact when $d(x, y)$ is
proper, since they are automatically bounded with respect to $d(x, y)$.
Conversely, if closed balls with respect to $d(x, y)$ are compact in $G$,
then $d(x, y)$ is proper on $G$, because closed subsets of compact sets
are compact.

        Suppose that $d(x, y)$ is a semi-ultrametric on a topological
group $G$ which is compatible with the topology on $G$, proper, and
invariant under left or right translations.  This implies that open
and closed balls in $G$ with respect to $d(x, y)$ centered at $e$ are
subgroups of $G$, as before.  Remember that open and closed balls with
respect to $d(x, y)$ are each both open and closed with respect to the
topology defined by $d(x, y)$, because $d(x, y)$ is a
semi-ultrametric.  Thus open and closed balls with respect to $d(x,
y)$ are each both open and closed with respect to the topology on $G$,
since $d(x, y)$ is compatible with the topology on $G$.  It follows
that open and closed balls with respect to $d(x, y)$ are compact
subsets of $G$, because $d(x, y)$ is also proper.  Every compact set
in $G$ is contained in a ball with respect to $d(x, y)$ centered at $e$,
and so $G$ has large compact open subgroups under these conditions.
Similarly, if $d(x, y)$ is invariant under both left and right translations, 
then $G$ has large compact open normal subgroups.

        Suppose now that $G$ is a topological group, and that
$U_1, U_2, U_3, \ldots$ is a sequence of compact open subgroups of $G$
such that $U_j \subseteq U_{j + 1}$ for each $j \in {\bf Z}_+$ and
$\bigcup_{j = 1}^\infty U_j = G$.  This implies that each compact set
$K \subseteq G$ is contained in the union of finitely many $U_j$'s,
and hence that $K \subseteq U_j$ for some $j$, so that $G$ has large
compact open subgroups.  If $U_j$ is also a normal subgroup of $G$ for
each $j$, then it follows that $G$ has large compact open normal
subgroups.  If $d(x, y)$ is a semi-ultrametric on $G$ which is
compatible with the topology on $G$, proper, and invariant under
left or right translations, as in the previous paragraph, then one
can take a sequence of balls centered at $e$ with respect to $d(x, y)$
with increasing radii tending to infinity.

        Let $r_j(x) = r_{U_j}(x)$ be the indicator function on $G$
associated to $U_j$, which is equal to $0$ when $x \in U_j$ and to $1$
when $x \in G$ is not in $U_j$.  Also let $t_1, t_2, t_3, \ldots$
be a monotone increasing sequence of nonnegative real numbers tending
to infinity, and put
\begin{equation}
\label{r(x) = max_{j ge 1} t_j r_j(x)}
        r(x) = \max_{j \ge 1} \, t_j \, r_j(x)
\end{equation}
for each $x \in G$.  By hypothesis, each $x \in G$ is contained in
$U_j$ for all but finitely many $j$, so that $r_j(x) = 0$ for all but
finitely many $j$, and hence the maximum in (\ref{r(x) = max_{j ge 1}
  t_j r_j(x)}) exists.  Clearly $r(x)$ satisfies $r(e) = 0$,
(\ref{r(x) = r(x^{-1})}), and (\ref{r(x y) le max(r(x), r(y))}),
because of the analogous properties of $r_j(x)$ for each $j$.  If
$U_j$ is a normal subgroup for each $j$, then $r_j(x)$ satisfies
(\ref{r(a x a^{-1}) = r(x)}) for each $j$, and hence $r(x)$ satisfies
(\ref{r(a x a^{-1}) = r(x)}) too.

        Thus $r(y^{-1} \, x)$ defines a left-invariant semi-ultrametric
on $G$ under these conditions, and $r(x \, y^{-1})$ defines a
right-invariant semi-ultrametric on $G$.  By construction, these
semi-ultrametrics are compatible with the topology on $G$.  These
semi-ultrametrics are also proper, because $U_j$ is compact for each
$j$ by hypothesis, and $t_j \to +\infty$ as $j \to \infty$.  If $U_j$
is a normal subgroup of $G$ for each $j$, then these two
semi-ultrametrics are the same, and hence invariant under both left
and right translations.

        Suppose that $G$ is a topological group with large compact
open subgroups which is also $\sigma$-compact, so that there is a
sequence $K_1, K_2, K_3, \ldots$ of compact subsets of $G$ such that
$\bigcup_{j = 1}^\infty K_j = G$.  Let $U_1$ be a compact open subgroup
of $G$ with $K_1 \subseteq U_1$, and for each $j \ge 2$, let $U_j$
be a compact open subgroup of $G$ such that $K_j \subseteq U_j$ and
$U_{j - 1} \subseteq U_j$.  This implies that $\bigcup_{j = 1}^\infty U_j
= G$, so that we are in the same situation as before.  If $G$ has large
compact open normal subgroups, then we can choose $U_j$ to be a normal
subgroup of $G$ for each $j$ too.

        Let $p$ be a prime number, and put
\begin{equation}
\label{d_p''(x, y) = |log_p |x|_p - log_p |y|_p|}
        d_p''(x, y) = \bigl|\log_p |x|_p - \log_p |y|_p\bigr|
\end{equation}
for every $x, y \in {\bf Q}_p^*$, where $\log_p t$ is the logarithm
base $p$ of a positive real number $t$.  It is easy to see that this
defines a translation-invariant semimetric on ${\bf Q}_p^*$ with
respect to multiplication, because $\log_p |x|_p$ is a homomorphism
from ${\bf Q}_p^*$ onto the additive group ${\bf Z}$ of integers.  One
can also check that (\ref{d_p''(x, y) = |log_p |x|_p - log_p |y|_p|})
is compatible with the usual topology on ${\bf Q}_p^*$ and proper.
However, (\ref{d_p''(x, y) = |log_p |x|_p - log_p |y|_p|}) is not a
semi-ultrametric on ${\bf Q}_p^*$, and indeed there is no
semi-ultrametric on ${\bf Q}_p^*$ with these properties, because ${\bf
  Q}_p^*$ does not have large compact subgroups.  Note that ${\bf Q}_p^*$
is isomorphic as a topological group with respect to multiplication
to the product of the additive group ${\bf Z}$ of integers and the
multiplicative group of $x \in {\bf Q}_p$ with $|x|_p = 1$.

\section{General linear groups}
\label{general linear groups}
\setcounter{equation}{0}

        Let $R$ be a commutative ring, and let $n$ be a positive integer.  
Also let $M_n(R)$ be the ring of $n \times n$ matrices with entries in
$R$, with respect to entry-wise addition and matrix multiplication.
Of course, $M_n(R)$ can be identified with $R^{n^2}$ in the usual way.
The determinant $\det A \in R$ of $A \in M_n(R)$ can be defined in the
standard way, and satisfies $\det (A \, B) = (\det A) \, (\det B)$ for
every $A, B \in M_n(R)$.

        Suppose that $R$ is a commutative topological ring.  This
means that $R$ is a commutative ring equipped with a topology which
makes it a topological group with respect to addition, and for which
multiplication is continuous as a mapping from $R \times R$ into $R$,
using the appropriate product topology on $R \times R$.  This implies
that $M_n(R)$ is a topological ring as well, with respect to the
topology that corresponds to the product topology on $R^{n^2}$.  Note
that the determinant is automatically continuous as a mapping from
$M_n(R)$ into $R$, since $\det A$ is a polynomial in the entries of
$A$.

        If $R$ has a nonzero multiplicative identity element $e$, then 
the identity matrix $I$ with diagonal entries equal to $e$ and
off-diagonal entries equal to $0$ is the multiplicative identity
element for $M_n(R)$, and $\det I = e$.  An element $A$ of $M_n(R)$ is
said to be invertible if it has a multiplicative inverse in $M_n(R)$,
in which case $\det A$ is invertible as an element of $R$.
Conversely, if $A \in M_n(R)$ and $\det A$ is invertible in $R$, then
$A$ is invertible in $M_n(R)$, as in Cramer's rule.  The group of
invertible elements of $M_n(R)$ is denoted $GL(n, R)$.

        If $R$ is a field, then $GL(n, R)$ consists of the $A \in M_n(R)$ 
such that $\det A \ne 0$.  Suppose that $R$ is a topological field, so
that $R$ is a topological ring, and $x \mapsto 1/x$ is continuous on
the set of $x \in R$ such that $x \ne 0$.  This implies that $GL(n,
R)$ is an open set in $M_n(R)$, and that the mapping from $A \in GL(n, R)$
to $A^{-1}$ is continuous, because of the expression for $A^{-1}$ in
terms of determinants.  It follows that $GL(n, R)$ is a topological
group with respect to the topology induced by the usual one on
$M_n(R)$, since matrix multiplication is continuous.  This includes
the cases where $R = {\bf R}$ or ${\bf C}$.

        Let $p$ be a prime number, so that the previous remarks can also
be applied to $R = {\bf Q}_p$.  Note that ${\bf Z}_p$ is a compact open
subring of ${\bf Q}_p$, and hence that $M_n({\bf Z}_p)$ is a compact
open subring of $M_n({\bf Q}_p)$.  An element $x$ of ${\bf Z}_p$ has
a multiplicative inverse in ${\bf Z}_p$ if and only if $|x|_p = 1$,
and the set of such $x$ is a compact open subset of ${\bf Q}_p$.
It follows that
\begin{equation}
\label{GL(n, {bf Z}_p) = {A in M_n({bf Z}_p) : |det A|_p = 1}}
        GL(n, {\bf Z}_p) = \{A \in M_n({\bf Z}_p) : |\det A|_p = 1\}
\end{equation}
is a compact open subgroup of $GL(n, {\bf Q}_p)$, since the
determinant is continuous as a mapping from $M_n({\bf Q}_p)$ into
${\bf Q}_p$, as before.

        Remember that $p^j \, {\bf Z}_p$ is a compact open subring of
${\bf Z}_p$ and hence ${\bf Q}_p$ for each positive integer $j$, so that
$M_n(p^j \, {\bf Z}_p)$ is a compact open subring of $M_n({\bf Z}_p)$ and
thus $M_n({\bf Q}_p)$.  Put
\begin{equation}
\label{GL_j(n, {bf Z}_p) = {A in M_n({bf Z}_p) : A - I in M_n(p^j {bf Z}_p)}}
 GL_j(n, {\bf Z}_p) = \{A \in M_n({\bf Z}_p) : A - I \in M_n(p^j \, {\bf Z}_p)\}
\end{equation}
for each $j \in {\bf Z}_+$.  If $A \in GL_j(n, {\bf Z}_p)$, then it is
easy to see that
\begin{equation}
\label{(det A) - 1 in p^j {bf Z}_p}
        (\det A) - 1 \in p^j \, {\bf Z}_p,
\end{equation}
because $\det I = 1$, and the determinant of $A$ is a sum of monomials
in the entries of $A$ with coefficients $\pm 1$.  In particular, this
implies that $|\det A|_p = 1$, so that $A \in GL(n, {\bf Z}_p)$.
One can also check that
\begin{equation}
\label{A^{-1} - I in M_n(p^j {bf Z}_p)}
        A^{-1} - I \in M_n(p^j \, {\bf Z}_p)
\end{equation}
under these conditions, which means that $A^{-1} \in GL_j(n, {\bf
  Z}_p)$.  The product of two elements of $GL_j(n, {\bf Z}_p)$ is in
$GL_j(n, {\bf Z}_p)$ as well, so that $GL_j(n, {\bf Z}_p)$ is a
subgroup of $GL(n, {\bf Z}_p)$.  More precisely, $GL_j(n, {\bf Z}_p)$
is a compact open subgroup of $GL(n, {\bf Z}_p)$ for each $j \in {\bf
  Z}_+$, which implies that $GL(n, {\bf Z}_p)$ and hence $GL(n, {\bf
  Q}_p)$ has small compact open subgroups.  In fact, $GL_j(n, {\bf
  Z}_p)$ is a normal subgroup of $GL(n, {\bf Z}_p)$ for each $j$, so
that $GL(n, {\bf Z}_p)$ has small compact open normal subgroups.

        Equivalently,
\begin{equation}
\label{M_n({bf Z}_p) = {A in M_n({bf Q}_p) : ||A|| le 1}}
        M_n({\bf Z}_p) = \{A \in M_n({\bf Q}_p) : \|A\| \le 1\},
\end{equation}
where $\|A\|$ is the ultranorm on $M_n({\bf Q}_p)$ defined in
(\ref{||A|| = max_{1 le j, k le n} |a_{j, k}|_p}), and
\begin{equation}
\label{GL(n, {bf Z}_p) = {A in GL(n, {bf Q}_p) : ||A||, ||A^{-1}|| le 1}}
        GL(n, {\bf Z}_p) = \{A \in GL(n, {\bf Q}_p) : \|A\|, \|A^{-1}\| \le 1\}.
\end{equation}
Alternatively,
\begin{equation}
\label{GL(n, {bf Z}_p) = {A in GL(n, {bf Q}_p) : ||A|| = ||A^{-1}|| = 1}}
        GL(n, {\bf Z}_p) = \{A \in GL(n, {\bf Q}_p) : \|A\| = \|A^{-1}\| = 1\},
\end{equation}
because
\begin{equation}
\label{1 = ||I|| le ||A|| ||A^{-1}||}
        1 = \|I\| \le \|A\| \, \|A^{-1}\|
\end{equation}
for every $A \in GL(n, {\bf Q}_p)$, by (\ref{||A B|| le ||A|| ||B||}).
Observe that
\begin{equation}
\label{||A C|| = ||C A|| = ||A||}
        \|A \, C\| = \|C \, A\| = \|A\|
\end{equation}
for every $A \in GL(n, {\bf Q}_p)$ and $C \in GL(n, {\bf Z}_p)$,
using (\ref{||A B|| le ||A|| ||B||}) again.  Thus
\begin{equation}
\label{||(A - B) C|| = ||C (A - B)|| = ||A - B||}
        \|(A - B) \, C\| = \|C \, (A - B)\| = \|A - B\|
\end{equation}
for every $A, B \in GL(n, {\bf Q}_p)$ and $C \in GL(n, {\bf Z}_p)$,
which implies in particular that the ultrametric $\|A - B\|$ is
invariant under both left and right translations on $GL(n, {\bf
  Z}_p)$.  Of course,
\begin{equation}
 GL_j(n, {\bf Z}_p) = \{A \in GL(n, {\bf Z}_p) : \|A - I\| \le p^{-j}\}
\end{equation}
for each $j \in {\bf Z}_+$.

        Put
\begin{equation}
\label{r(A) = max(||A - I||, ||A^{-1} - I||)}
        r(A) = \max(\|A - I\|, \|A^{-1} - I\|)
\end{equation}
for every $A \in GL(n, {\bf Q}_p)$.  If $A \in GL(n, {\bf Z}_p)$, then
\begin{equation}
\label{r(A) = ||A - I|| = ||A^{-1} - I|| le 1}
        r(A) = \|A - I\| = \|A^{-1} - I\| \le 1.
\end{equation}
Otherwise, if $A \in GL(n, {\bf Q}_p)$ is not in $GL(n, {\bf Z}_p)$,
then either $\|A\| > 1$ or $\|A^{-1}\| > 1$, and hence
\begin{equation}
\label{r(A) = max(||A||, ||A^{-1}||) ge p}
        r(A) = \max(\|A\|, \|A^{-1}\|) \ge p,
\end{equation}
because $\|\cdot \|$ is an ultranorm on $M_n({\bf Q}_p)$.  If $A, B
\in GL(n, {\bf Z}_p)$, then
\begin{eqnarray}
\label{r(A B) = ||A B - I|| le ... max(||A - I||, ||B - I||) = max(r(A), r(B))}
 r(A \, B) & = & \|A \, B - I\| \le \max(\|A \, B - B\|, \|B - I\|) \\
 & = & \max(\|A - I\|, \|B - I\|) = \max(r(A), r(B)). \nonumber
\end{eqnarray}
using (\ref{||(A - B) C|| = ||C (A - B)|| = ||A - B||}) in the third
step.  Note that $r(A^{-1}) = r(A)$ for each $A \in GL(n, {\bf Q}_p)$,
$r(I) = 0$, and that $r(A) > 0$ when $A \ne I$.

        Let $t$ be a real number such that $1 \le t \le p$, and put
\begin{equation}
\label{r'(A) = min(r(A), t)}
        r'(A) = \min(r(A), t)
\end{equation}
for every $A \in GL(n, {\bf Q}_p)$.  Thus $r'(A) = r(A)$ when $A \in
GL(n, {\bf Z}_p)$, and $r'(A) = t$ otherwise.  It is easy to see that
\begin{equation}
\label{r'(A B) le max(r'(A), r'(B))}
        r'(A \, B) \le \max(r'(A), r'(B))
\end{equation}
for every $A, B \in GL(n, {\bf Q}_p)$, using (\ref{r(A B) = ||A B -
  I|| le ... max(||A - I||, ||B - I||) = max(r(A), r(B))}) when $A, B
\in GL(n, {\bf Z}_p)$.  As before, $r'(A^{-1}) = r(A)$ for every
$A \in GL(n, {\bf Q}_p)$, $r'(I) = 0$, and $r'(A) > 0$ when $A \ne I$.

        As in Section \ref{topological groups}, $r'(B^{-1} \, A)$ defines
a left-invariant ultrametric on $GL(n, {\bf Q}_p)$, and $r'(A \, B^{-1})$
defines a right-invariant ultrametric on $GL(n, {\bf Q}_p)$.  If $A, B
\in GL(n, {\bf Z}_p)$, then
\begin{equation}
\label{r'(B^{-1} A) = r(B^{-1} A) = ||B^{-1} A - I|| = ||A - B||}
        r'(B^{-1} \, A) = r(B^{-1} \, A) = \|B^{-1} \, A - I\| = \|A - B\|,
\end{equation}
and similarly for $r'(A \, B^{-1})$.  One can check that the
topologies on $GL(n, {\bf Q}_p)$ determined by these metrics are the
same as the standard topology, using translation-invariance and direct
comparisons around the identity matrix $I$.

        Put
\begin{equation}
\label{r''(A) = max(log_p ||A||, log_p ||A^{-1}||)}
        r''(A) = \max(\log_p \|A\|, \log_p \|A^{-1}\|)
\end{equation}
for each $A \in GL(n, {\bf Q}_p)$.  Note that $\|A\|$ is an integer
power of $p$ for every $A \in M_n({\bf Q}_p)$ with $A \ne 0$, so that
$\log_p \|A\| \in {\bf Z}$.  If $A \in GL(n, {\bf Q}_p)$, then $\|A\|
\ge 1$ or $\|A^{-1}\| \ge 1$, by (\ref{1 = ||I|| le ||A||
  ||A^{-1}||}), which implies that $r''(A) \ge 0$.  By construction,
$r''(A^{-1}) = r''(A)$ for every $A \in GL(n, {\bf Q}_p)$, and $r''(A)
= 0$ if and only if $A \in GL(n, {\bf Z}_p)$.  It is easy to see that
\begin{equation}
\label{r''(A B) le r''(A) + r''(B)}
        r''(A \, B) \le r''(A) + r''(B)
\end{equation}
for every $A, B \in GL(n, {\bf Q}_p)$, using (\ref{||A B|| le ||A||
  ||B||}).  Thus $r''(B^{-1} \, A)$ defines a left-invariant
semimetric on $GL(n, {\bf Q}_p)$, and $r''(A \, B^{-1})$ defines a
right-invariant semimetric on $GL(n, {\bf Q}_p)$, as in Section
\ref{topological groups}.  One can also check that these semimetrics
are compatible with the standard topology on $GL(n, {\bf Q}_p)$,
and that they are proper.  Of course, $GL(n, {\bf Q}_p)$ does not
have large compact subgroups, even when $n = 1$.

\section{Heisenberg groups}
\label{heisenberg groups}
\setcounter{equation}{0}

        Let $R$ be a commutative ring, let $n$ be a positive integer, and 
put
\begin{equation}
\label{H_n(R) = R^n times R^n times R}
        H_n(R) = R^n \times R^n \times R,
\end{equation}
initially as a set.  If $(x, y, t), (x', y', t') \in H_n(R)$, so that
$x, y, x', y' \in R^n$ and $t, t' \in R$, then put
\begin{equation}
\label{(x, y, t) diamond (x', y', t') = ...}
        (x, y, t) \diamond (x', y', t') = 
           \Big(x + x', y + y', t + t' + \sum_{j = 1}^n x_j \, y_j'\Big).
\end{equation}
If $(x'', y'', t'') \in H_n(R)$ too, then
\begin{eqnarray}
\label{((x, y, t) diamond (x', y', t')) diamond (x'', y'', t'') = ...}
\lefteqn{\quad  \qquad  ((x, y, t) \diamond (x', y', t')) 
                                              \diamond (x'', y'', t'')} \\
 & = & \Big(x + x', y + y', t + t' + \sum_{j = 1}^n x_j \, y_j'\Big)
                                    \diamond (x'', y'', t'') \nonumber \\
 & = & \Big(x + x' + x'', y + y' + y'', t + t' + t'' + \sum_{j = 1}^n x_j \, y_j'
  + \sum_{j = 1}^n x_j \, y_j'' + \sum_{j = 1}^n x_j' \, y_j''\Big). \nonumber
\end{eqnarray}
Similarly,
\begin{eqnarray}
\label{(x, y, t) diamond ((x', y', t') diamond (x'', y'', t''))}
\lefteqn{\quad  \qquad  (x, y, t) \diamond ((x', y', t')
                                                \diamond (x'', y'', t''))} \\
 & = & (x, y, t) \diamond \Big(x' + x'', y' + y'', t' + t''
                          + \sum_{j = 1}^n x_j' \, y_j''\Big) \nonumber \\
 & = & \Big(x + x' + x'', y + y' + y'', t + t' + t'' + \sum_{j = 1}^n x_j \, y_j'
 + \sum_{j = 1}^n x_j \, y_j'' + \sum_{j = 1}^n x_j' \, y_j''\Big). \nonumber
\end{eqnarray}
Thus (\ref{((x, y, t) diamond (x', y', t')) diamond (x'', y'', t'') =
  ...}) and (\ref{(x, y, t) diamond ((x', y', t') diamond (x'', y'',
  t''))}) are equal to each other, which implies that $\diamond$ is
associative on $H_n(R)$.  Observe that
\begin{equation}
\label{(0, 0, 0) diamond (x, y, t) = (x, y, t) diamond (0, 0, 0) = (x, y, t)}
 (0, 0, 0) \diamond (x, y, t) = (x, y, t) \diamond (0, 0, 0) = (x, y, t)
\end{equation}
for every $(x, y, t) \in H_n(R)$, so that $(0, 0, 0)$ is the identity element
in $H_n(R)$ with respect to $\diamond$.  If we put
\begin{equation}
\label{(x, y, t)^{-1} = (-x, -y, -t + sum_{j = 1}^n x_j y_j)}
        (x, y, t)^{-1} = \Big(-x, -y, -t + \sum_{j = 1}^n x_j \, y_j\Big)
\end{equation}
for every $(x, y, t) \in H_n(R)$, then it is easy to see that
\begin{equation}
\label{(x, y, t) diamond (x, y, t)^{-1} = ... = (0, 0, 0)}
        (x, y, t) \diamond (x, y, t)^{-1} = (x, y, t)^{-1} \diamond (x, y, t) 
                                          = (0, 0, 0),
\end{equation}
so that $H_n(R)$ is a group with respect to $\diamond$.  If $R$ is a
commutative topological ring, then $H_n(R)$ is a topological group,
with respect to the corresponding product topology.

        If $r \in R$ and $(x, y, t) \in H_n(R)$, then put
\begin{equation}
\label{delta_r((x, y, t)) = (r x, r y, r^2 t)}
        \delta_r((x, y, t)) = (r \, x, r \, y, r^2 \, t).
\end{equation}
Thus
\begin{eqnarray}
\label{delta_r((x, y, t)) diamond delta_r((x', y', t')) ...}
\lefteqn{\delta_r((x, y, t)) \diamond \delta_r((x', y', t'))} \\
 & = & (r \, x, r \, y, r^2 \, t) \diamond (r \, x', r \, y', r^2 \, t') 
                                                              \nonumber \\
 & = & \Big(r \, x + r \, x', r \, y + r \, y', r^2 \, t + r^2 \, t'
                              + r^2 \sum_{j = 1}^n x_j \, y_j'\Big) \nonumber \\
 & = & \delta_r((x, y, t) \diamond (x', y', t'))   \nonumber
\end{eqnarray}
for every $(x, y, t), (x', y', t') \in H_n(R)$ and $r \in R$, so that
$\delta_r$ defines a group homomorphism from $H_n(R)$ into itself.
This implies that $\delta_r((0, 0, 0)) = (0, 0, 0)$ for every $r \in
R$, and that
\begin{equation}
\label{delta_r((x, y, t)^{-1}) = (delta_r((x, y, t)))^{-1}}
        \delta_r((x, y, t)^{-1}) = (\delta_r((x, y, t)))^{-1}
\end{equation}
for every $(x, y, t) \in H_n(R)$ and $r \in R$, which can also be
verified directly.  Note that
\begin{equation}
\label{delta_r circ delta_{r'} = delta_{r r'}}
        \delta_r \circ \delta_{r'} = \delta_{r \, r'}
\end{equation}
for every $r, r' \in R$, and that $\delta_r((x, y, t)) = (0, 0, 0)$
for every $(x, y, t) \in H_n(R)$ when $r = 0$.  If $R$ has a
multiplicative identity element $e$, then $\delta_e$ is the identity
mapping on $H_n(R)$.  In this case, if $r \in R$ has a multiplicative
inverse, then $\delta_r$ is an invertible mapping on $H_n(R)$, with
inverse equal to $\delta_{r^{-1}}$.  In particular, if $R$ is a field,
then $\delta_r$ is a one-to-one mapping from $H_n(R)$ onto itself for
every $r \in R$ with $r \ne 0$.  If $R$ is a commutative topological
ring, then $\delta_r$ is a continuous mapping from $H_n(R)$ into
itself for each $r \in R$, and indeed $\delta_r((x, y, t))$ is
continuous as a function of $r \in R$ and $(x, y, t) \in H_n(R)$.  It
follows that $\delta_r$ is a homeomorphism from $H_n(R)$ onto itself
when $r \in R$ has a multiplicative inverse in $R$.

        Let $p$ be a prime number, and let us apply the previous discussion
to $R = {\bf Q}_p$.  We can also apply the previous discussion to 
$R = {\bf Z}_p$, to get a group $H_n({\bf Z}_p)$ which is a compact
open subgroup of $H_n({\bf Q}_p)$.  Similarly, if $k$, $l$ are integers such 
that $2 \, k \ge l$, then it is easy to see that
\begin{equation}
\label{(p^k {bf Z}_p)^n times (p^k {bf Z}_p)^n times (p^l {bf Z}_p)}
 (p^k \, {\bf Z}_p)^n \times (p^k \, {\bf Z}_p)^n \times (p^l \, {\bf Z}_p)
\end{equation}
is a compact open subgroup of $H_n({\bf Q}_p)$.  It follows that
$H_n({\bf Q}_p)$ has large and small compact open subgroups.

        If $R$ is any commutative ring and $(x, y, t), (x', y', t') 
\in H_n(R)$, then
\begin{eqnarray}
\label{((x, y, t) diamond (x', y', t')) diamond (x, y, t)^{-1} = ...}
\lefteqn{\quad  ((x, y, t) \diamond (x', y', t')) \diamond (x, y, t)^{-1}} \\
 & = & \Big(x + x', y + y', t + t' + \sum_{j = 1}^n x_j \, y_j'\Big)
        \diamond \Big(-x, -y, -t + \sum_{j = 1}^n x_j \, y_j\Big) \nonumber \\
 & = & \Big(x', y', t' + \sum_{j = 1}^n x_j \, y_j' + \sum_{j = 1}^n x_j \, y_j
                         - \sum_{j = 1}^n (x_j + x_j') \, y_j\Big) \nonumber \\
 & = & \Big(x', y', t' + \sum_{j = 1}^n (x_j \, y_j' - x_j' \, y_j)\Big).
                                                                  \nonumber
\end{eqnarray}
Observe that
\begin{equation}
\label{{0} times {0} times R}
        \{0\} \times \{0\} \times R
\end{equation}
is a normal subgroup of $H_n(R)$, and that every element of (\ref{{0}
  times {0} times R}) commutes with every element of $H_n(R)$.
Suppose that $R$ is a field, and that $K$ is a normal subgroup of
$H_n(R)$.  If $(x', y', t') \in K$ and $x' \ne 0$ or $y' \ne 0$, then
it follows from (\ref{((x, y, t) diamond (x', y', t')) diamond (x, y,
  t)^{-1} = ...}) that $(x', y', t'') \in K$ for every $t '' \in R$.
This implies that (\ref{{0} times {0} times R}) is contained in $K$,
by multiplying by $(x', y', t')^{-1}$.

         Applying this to $R = {\bf Q}_p$ again, we get that $H_n({\bf Q}_p)$
has no nontrivial compact normal subgroups, and that $H_n({\bf Q}_p)$
does not have small open normal subgroups.  However,
\begin{equation}
\label{H_n(p^k {bf Z}_p) = ...}
        H_n(p^k \, {\bf Z}_p) = 
    (p^k \, {\bf Z}_p)^n \times (p^k \, {\bf Z}_p)^n \times (p^k \, {\bf Z}_p)
\end{equation}
is a compact open normal subgroup of $H_n({\bf Z}_p)$ for each
nonnegative integer $k$, by (\ref{((x, y, t) diamond (x', y', t'))
  diamond (x, y, t)^{-1} = ...}).  In particular, $H_n({\bf Z}_p)$ has
small compact open normal subgroups.  Note that $p^k \, {\bf Z}_p$ is
an ideal in ${\bf Z}_p$ and a sub-ring of ${\bf Q}_p$ when $k \ge 0$,
and that (\ref{H_n(p^k {bf Z}_p) = ...}) is the same as (\ref{(p^k {bf
    Z}_p)^n times (p^k {bf Z}_p)^n times (p^l {bf Z}_p)}) with $l =
k$, which satisfies $2 \, k \ge l$ when $k \ge 0$.

        Put
\begin{equation}
\label{N((x, y, t)) = ...}
        N((x, y, t)) = 
           \max(|x_1|_p, \ldots, |x_n|_p, |y_1|, \ldots, |y_n|_p, |t|_p^{1/2})
\end{equation}
for each $(x, y, t) \in H_n({\bf Q}_p)$.  One can check that
\begin{equation}
\label{N((x, y, t) diamond (x', y', t')) le ...}
 N((x, y, t) \diamond (x', y', t')) \le \max(N((x, y, t)), N((x', y', t')))
\end{equation}
for every $(x, y, t), (x', y', t') \in H_n({\bf Q}_p)$, using the
simple fact that
\begin{equation}
\label{(a b)^{1/2} le max(a, b)}
        (a \, b)^{1/2} \le \max(a, b)
\end{equation}
for any two nonnegative real numbers $a$, $b$.  Similarly,
\begin{equation}
\label{N((x, y, t)^{-1}) = N((x, y, t))}
        N((x, y, t)^{-1}) = N((x, y, t))
\end{equation}
for every $(x, y, t) \in H_n({\bf Q}_p)$.  By construction,
\begin{equation}
\label{N(delta_r((x, y, t))) = |r|_p N((x, y, t))}
        N(\delta_r((x, y, t))) = |r|_p \, N((x, y, t))
\end{equation}
for every $(x, y, t) \in H_n({\bf Q}_p)$ and $r \in {\bf Q}_p$.
It follows that
\begin{equation}
\label{N((x', y', t')^{-1} diamond (x, y, t))}
        N((x', y', t')^{-1} \diamond (x, y, t))
\end{equation}
defines a left-invariant unltrametric on $H_n({\bf Q}_p)$, and that
\begin{equation}
\label{N((x, y, t) diamond (x', y', t')^{-1})}
        N((x, y, t) \diamond (x', y', t')^{-1})
\end{equation}
defines a right-invariant ultrametric on $H_n({\bf Q}_p)$, as in
Section \ref{topological groups}.  Of course, the topologies on
$H_n({\bf Q}_p)$ determined by these ultrametrics are the same as the
product topology associated to the standard topology on ${\bf Q}_p$.
These ultrametrics also transform in a nice way under the dilations
$\delta_r$ for each $r \in {\bf Q}_p$, by (\ref{N(delta_r((x, y, t)))
  = |r|_p N((x, y, t))}).  Observe that
\begin{eqnarray}
\label{{(x, y, t) in H_n({bf Q}_p) : N((x, y, t)) le p^{-k}} = ...}
\lefteqn{\{(x, y, t) \in H_n({\bf Q}_p) : N((x, y, t)) \le p^{-k}\}} \\
 & = & (p^k \, {\bf Z}_p)^n \times (p^k \, {\bf Z}_p)^n 
                             \times (p^{2 k} \, {\bf Z}_p) \nonumber \\
 & = & \delta_{p^k}(H_n({\bf Z}_p))                         \nonumber
\end{eqnarray}
for each $k \in {\bf Z}$, which is the same as (\ref{(p^k {bf Z}_p)^n
  times (p^k {bf Z}_p)^n times (p^l {bf Z}_p)}) with $l = 2 \, k$.

        Alternatively, put
\begin{equation}
\label{widetilde{N}((x, y, t)) = ...}
        \widetilde{N}((x, y, t)) =
         \max(|x_1|_p, \ldots, |x_n|_p, |y_1|_p, \ldots, |y_n|_p, |t|_p)
\end{equation}
for every $(x, y, t) \in H_n({\bf Z}_p)$.  As before,
\begin{equation}
\label{widetilde{N}((x, y, t) diamond (x', y', t')) le ...}
        \widetilde{N}((x, y, t) \diamond (x', y', t'))
         \le \max(\widetilde{N}((x, y, t)), \widetilde{N}((x', y', t')))
\end{equation}
for every $(x, y, t), (x', y', t') \in H_n({\bf Z}_p)$, and
\begin{equation}
\label{widetilde{N}((x, y, t)^{-1}) = widetilde{N}((x, y, t))}
        \widetilde{N}((x, y, t)^{-1}) = \widetilde{N}((x, y, t))
\end{equation}
for every $(x, y, t) \in H_n({\bf Z}_p)$.  One can also check that
\begin{equation}
\label{widetilde{N}(((x, y, t) diamond (x', y', t')) diamond (x, y, t)^{-1})}
 \widetilde{N}(((x, y, t) \diamond (x', y', t')) \diamond (x, y, t)^{-1})
        = \widetilde{N}((x', y', t'))
\end{equation}
for every $(x, y, t), (x', y', t') \in H_n({\bf Z}_p)$, using
(\ref{((x, y, t) diamond (x', y', t')) diamond (x, y, t)^{-1} = ...}).
This implies that
\begin{equation}
\label{widetilde{N}((x', y', t')^{-1} diamond (x, y, t)) = ...}
        \widetilde{N}((x', y', t')^{-1} \diamond (x, y, t))
           = \widetilde{N}((x, y, t) \diamond (x', y', t')^{-1})
\end{equation}
defines an ultrametric on $H_n({\bf Z}_p)$ which is invariant under
both left and right translations, as in Section \ref{topological
  groups}.  As usual, the topology on $H_n({\bf Z}_p)$ determined by
this ultrametric is the same as the product topology corresponding to
the standard topology on ${\bf Z}_p$.  Note that
\begin{equation}
\label{{(x, y, t) in H_n({bf Z}_p) : widetilde{N}((x, y, t)) le p^{-k}} = ...}
        \{(x, y, t) \in H_n({\bf Z}_p) : \widetilde{N}((x, y, t)) \le p^{-k}\}
                                                        = H_n(p^k \, {\bf Z}_p)
\end{equation}
for each nonngative integer $k$, which is the same as (\ref{H_n(p^k
  {bf Z}_p) = ...}).  Observe also that
\begin{equation}
\label{N((x, y, t))^2 le widetilde{N}((x, y, t)) le N((x, y, t))}
        N((x, y, t))^2 \le \widetilde{N}((x, y, t)) \le N((x, y, t))
\end{equation}
for every $(x, y, t) \in H_n({\bf Z}_p)$, which reflects the fact that
the ultrametrics (\ref{N((x', y', t')^{-1} diamond (x, y, t))}),
(\ref{N((x, y, t) diamond (x', y', t')^{-1})}), and
(\ref{widetilde{N}((x', y', t')^{-1} diamond (x, y, t)) = ...})
determine the same topology on $H_n({\bf Z}_p)$.

\section{Upper-triangular matrices}
\label{upper-triangular matrices}
\setcounter{equation}{0}

        Let $R$ be a commutative ring, let $n$ be a positive integer,
and let $M_n(R)$ be the ring of $n \times n$ matrices $A = \{a_{j,
  k}\}_{j, k = 1}^n$ with entries in $R$, as usual.  More precisely,
the index $j$ indicates the row in the matrix, and the index $k$
indicates the column in the matrix, so that $A$ looks as follows.
\begin{equation}
\label{matrix in M_n(R)}
\left(\begin{array}{ccccc}
        a_{1, 1}      & a_{1, 2}    & \cdots & a_{1, n - 1}    & a_{1, n} \\
        a_{2, 1}      & a_{2, 2}    & \cdots & a_{2, n - 2}    & a_{2, n} \\
        \vdots       & \vdots     & \ddots & \vdots         & \vdots    \\
        a_{n - 1, 1}  & a_{n - 1, 2} & \cdots & a_{n - 1, n - 1} & a_{n - 1, n} \\
        a_{n, 1}      & a_{n, 2}    & \cdots & a_{n, n - 1}    & a_{n, n} \\
                                                     \end{array}\right)
\end{equation}
Let $T_n(R)$ be the collection of upper-triangular matrices $A$, which
means that $a_{j, k} = 0$ when $j > k$, as below.
\begin{equation}
\label{matrix in T_n(R)}
\left(\begin{array}{ccccc}
        a_{1, 1} & a_{1, 2} & \cdots & a_{1, n - 1}    & a_{1, n} \\
        0       & a_{2, 2} & \cdots & a_{2, n - 1}    & a_{2, n} \\
        \vdots  & \vdots  & \ddots & \vdots         & \vdots  \\
        0       & 0       & \cdots & a_{n - 1, n - 1} & a_{n - 1, n} \\
        0       & 0       & \cdots & 0              & a_{n, n} \\
                                                 \end{array}\right)
\end{equation}
It is well known and easy to see that $T_n(A)$ is a subring of $M_n(R)$,
and that the determinant of (\ref{matrix in T_n(R)}) is equal to the
product of its diagonal entries.

        Let $l$ be a nonnegative integer, and put
\begin{equation}
\label{T_n^l(R) = ...}
        T_n^l(R) = \{A = \{a_{j, k}\}_{j, k = 1}^n \in T_n(R) : 
                               a_{j, k} = 0 \hbox{ when } k - j \ne l\}.
\end{equation}
Thus $T_n^0(R)$ consists of diagonal matrices, $T_n^l(R) = \{0\}$ when
$l \ge n$, and $T_n^l(R)$ is a subgroup of $T_n(R)$ as a commutative
group with respect to addition for each $l$.  In fact, $T_n(R)$ is the
direct sum of $T_n^l(R)$ for $l = 0, 1, \ldots, n - 1$ as a
commutative group with respect to addition.  If $A \in T_n^l(R)$ and
$A' \in T_n^{l'}(R)$ for some nonnegative integers $l$, $l'$, then
\begin{equation}
\label{A A' in T_n^{l + l'}(R)}
        A \, A' \in T_n^{l + l'}(R),
\end{equation}
and in particular $A \, A' = 0$ when $l + l' \ge n$.

        Let $r \in R$ and $A = \{a_{j, k}\}_{j, k = 1}^n \in T_n(R)$ be given,
and let $\delta_r(A)$ be the element of $T_n(R)$ whose $(j, k)$th entry is
equal to
\begin{equation}
\label{r^{k - j} a_{j, k}}
        r^{k - j} \, a_{j, k}
\end{equation}
when $j < k$, and to $a_{j, j}$ when $j = k$.  Of course, the $(j,
k)$th entry of $\delta_r(A)$ is $0$ when $j > k$.  If $r = 0$, then
$\delta_r(A)$ is the diagonal matrix whose diagonal entries are the
same as the diagonal entries of $A$.  If $R$ has a nonzero
multiplicative identity element and $r$ has a multiplicative inverse,
then $\delta_r(A)$ can be obtained by conjugating $A$ by a diagonal
matrix whose entries are given by successive powers of $r$.  Of
course,
\begin{equation}
\label{delta_r(A + A') = delta_r(A) + delta_r(A')}
        \delta_r(A + A') = \delta_r(A) + \delta_r(A')
\end{equation}
for every $A, A' \in T_n(R)$ and $r \in R$.

        Equivalently, if $r$ is any element of $R$ and $A \in T_n^l(R)$
for some nonnegative integer $l$, then
\begin{equation}
\label{delta_r(A) = r^l A}
        \delta_r(A) = r^l \, A
\end{equation}
when $l > 0$, and $\delta_r(A) = A$ when $l = 0$.  Using this, it is easy
to see that
\begin{equation}
\label{delta_r(A A') = delta_r(A) delta_r(A')}
        \delta_r(A \, A') = \delta_r(A) \, \delta_r(A')
\end{equation}
for every $A, A' \in T_n(R)$, by reducing to the case where $A \in
T_n^l(R)$, $A' \in T_n^{l'}(R)$ for some $l, l' \ge 0$.
Alternatively, if $R$ has a multiplicative identity element and $r$
has a multiplicative inverse in $R$, then (\ref{delta_r(A A') =
  delta_r(A) delta_r(A')}) follows by expressing $\delta_r$ in terms
of conjugation, as in the previous paragraph.  At any rate, $\delta_r$
is a ring homomorphism from $T_n(R)$ into itself for each $r \in R$.

        Clearly
\begin{equation}
\label{delta_r circ delta_{r'} = delta_{r r'}, 2}
        \delta_r \circ \delta_{r'} = \delta_{r \, r'}
\end{equation}
for every $r, r' \in R$.  If $R$ has a nonzero multiplicative identity
element $e$, then $\delta_e$ is the identity mapping on $T_n(R)$.  In
this case, if $r \in R$ has a multiplicative inverse in $R$, then
$\delta_r$ is an invertible mapping on $T_n(R)$, with inverse
$\delta_{r^{-1}}$.

        If $R$ is a topological ring, then $T_n(R)$ is a closed set in
$M_n(R)$, with respect to the usual topology on $M_n(R)$.  Of course,
we can identify $T_n(R)$ with the Cartesian product of $n (n + 1) / 2$
copies of $R$, and the topology on $T_n(R)$ induced by the usual
topology on $M_n(R)$ is the same as the product topology on this
Cartesian product associated to the given topology on $R$.  In this
case, $\delta_r$ is a continuous mapping from $T_n(A)$ into itself for
each $r \in R$, and in fact $\delta_r(A)$ is continuous as a function
of $r \in R$ and $A \in T_n(A)$.  In particular, $\delta_r$ is a
homeomorphism from $T_n(A)$ onto itself when $R$ has a multiplicative
identity element and $r \in R$ has a multiplicative inverse in $R$.

        Let $T_n^+(R)$ be the collection of matrices $A = 
\{a_{j, k}\}_{j, k = 1}^n \in M_n(R)$ which are strictly upper-triangular,
in the sense that $a_{j, k} = 0$ when $j \ge k$.  This is is a subring
of $T_n(R)$, which is the same as the direct sum of $T_n^l(R)$ for $l
= 1, \ldots, n - 1$ as a commutative group with respect to addition.  
Note that $\delta_0(A) = 0$ and $A^n = 0$ for every $A \in T_n^+(R)$.
If $R$ is a topological ring, then $T_n^+(R)$ is also a closed set in
$M_n(R)$, which can be identified with the Cartesian product of
$n (n - 1) / 2$ copies of $R$ with the corresponding product topology.

        Now let $p$ be a prime number, and let us apply the previous
remarks to $R = {\bf Q}_p$.  Put
\begin{equation}
\label{N(A) = max_{1 le j < k le n} |a_{j, k}|_p^{1/(k - j)}}
        N(A) = \max_{1 \le j < k \le n} |a_{j, k}|_p^{1/(k - j)}
\end{equation}
for each $A = \{a_{j, k}\}_{j, k = 1}^n \in T_n({\bf Q}_p)$, so that
$N(A) = 0$ if and only if $A \in T_n^0({\bf Q}_p)$.  It is easy to see that
\begin{equation}
\label{N(A + A') le max(N(A), N(A'))}
        N(A + A') \le \max(N(A), N(A'))
\end{equation}
for every $A, A' \in T_n({\bf Q}_p)$, because of the ultrametric
version of the triangle inequality for the $p$-adic absolute value.
By construction,
\begin{equation}
\label{N(delta_r(A)) = |r|_p N(A)}
        N(\delta_r(A)) = |r|_p \, N(A)
\end{equation}
for every $A \in T_n({\bf Q}_p)$ and $r \in {\bf Q}_p$.

        Suppose that $A = \{a_{j, k}\}_{j, k = 1}^n, A' = \{a_{j, k}\}_{j, k = 1}^n
\in T_n({\bf Q}_p)$ have their diagonal entries in ${\bf Z}_p$, and let us 
check that
\begin{equation}
\label{N(A A') le max(N(A), N(A'))}
        N(A \, A') \le \max(N(A), N(A')).
\end{equation}
It suffices to show that
\begin{equation}
\label{|a_{j, k} a'_{k, l}|_p^{1/(l - j)} le max(N(A), N(A'))}
        |a_{j, k} \, a'_{k, l}|_p^{1/(l - j)} \le \max(N(A), N(A'))
\end{equation}
when $1 \le j < l \le n$ and $j \le k \le l$, by the definition of
matrix multiplication, and the ultrametric version of the triangle
inequality for the $p$-adic absolute value.  If $j = k$, then
\begin{equation}
\label{|a_{j, j} a'_{j, l}|_p^{1/(l - j)} le |a'_{j, l}|_p^{1/(l - j)} le N(A')}
 |a_{j, j} \, a'_{j, l}|_p^{1/(l - j)} \le |a'_{j, l}|_p^{1/(l - j)} \le N(A'),
\end{equation}
as desired, because $a_{j, j} \in {\bf Z}_p$ by hypothesis.  Similarly,
if $k = l$, then
\begin{equation}
\label{|a_{j, l} a'_{l, l}|_p^{1/(l - j)} le |a_{j, l}|_p^{1/(l - j)} le N(A)}
 |a_{j, l} \, a'_{l, l}|_p^{1/(l - j)} \le |a_{j, l}|_p^{1/(l - j)} \le N(A),
\end{equation}
as desired, because $a'_{l, l} \in {\bf Z}_p$ by hypothesis.  
Otherwise, if $j < k < l $, then
\begin{equation}
\label{|a_{j, k} a'_{k, l}|_p le ... le max(N(A), N(A'))^{l - j}}
        |a_{j, k} \, a'_{k, l}|_p \le N(A)^{k - j} \, N(A')^{l - k}
                                 \le \max(N(A), N(A'))^{l - j},
\end{equation}
which implies (\ref{|a_{j, k} a'_{k, l}|_p^{1/(l - j)} le max(N(A),
  N(A'))}) also in this case.

\section{Upper-triangular matrices, 2}
\label{upper-triangular matrices, 2}
\setcounter{equation}{0}

        Let $R$ be a commutative ring with nonzero multiplicative 
identity element $e$, and let $n$ be a positive integer.  Also let
$M_n(R)$ be the ring of $n \times n$ matrices $A = \{a_{j, k}\}_{j, k
  = 1}^n$ with entries in $R$, and let $T_n(R)$ and $T_n^+(R)$ be the
sub-rings of $M_n(R)$ consisting of upper-triangular and strictly
upper-triangular matrices, respectively, as in the previous section.
Consider the collection $T^+(n, R)$ of $A = \{a_{j, k}\}_{j, k = 1}^n
\in T_n(R)$ whose diagonal entries are equal to $e$, so that $A$ has
the following form.
\begin{equation}
\label{matrix in T^+(n, R)}
\left(\begin{array}{ccccc} e      & a_{1, 2} & \cdots & a_{1, n - 1} & a_{1, n} \\
                           0      & e       & \cdots & a_{2, n - 1} & a_{2, n} \\
                           \vdots & \vdots  & \ddots & \vdots      & \vdots \\
                           0      & 0       & \cdots & e       & a_{n - 1, n} \\
                           0      & 0       & \cdots & 0       & e          \\
                                                           \end{array}\right)
\end{equation}
Equivalently, $a_{j, j} = e$ for $j = 1, \ldots, n$ and $a_{j, k} = 0$
when $j > k$, which is the same as saying that $A - I \in T_n^+(R)$.
Thus $\det A = e$ for every $A \in T^+(n, R)$, and it is well known
that $T^+(n, R)$ is a subgroup of $GL(n, R)$.  If $R$ is a topological
ring, then $T^+(n, R)$ is a closed set in $M_n(R)$, with respect to
the usual topology on $M_n(R)$.  In particular, $T^+(n, R)$ is a
relatively closed subgroup of $GL(n, R)$, with respect to the topology
induced by the usual topology on $M_n(R)$.

        Of course, $T^+(1, R)$ is the trivial group, consisting of the 
$1 \times 1$ matrix whose only entry is equal to $e$.  Every element
of $T^+(2, R)$ is of the form
\begin{equation}
\label{R into T^+(2, R)}
\left(\begin{array}{cc} e & a \\ 0 & e \end{array}\right)
\end{equation}
for some $a \in R$, and
\begin{equation}
\label{multiplication in T^+(2, R)}
 \left(\begin{array}{cc} e & a \\ 0 & e \end{array}\right)
  \left(\begin{array}{cc} e & b \\ 0 & e \end{array}\right)
   = \left(\begin{array}{cc} e & a + b \\ 0 & e \end{array}\right)
\end{equation}
for every $a, b \in R$.  Thus the mapping from $a \in R$ to (\ref{R
  into T^+(2, R)}) defines a group isomorphism from $R$ as a
commutative group with respect to addition onto $T^+(2, R)$ as a group
with respect to matrix multiplication.  Similarly, every element of
$T^+(3, R)$ can be expressed as
\begin{equation}
\label{H_1(R) into T^+(3, R)}
\left(\begin{array}{ccc} e & x & t \\ 0 & e & y \\ 0 & 0 & e \end{array}\right)
\end{equation}
for some $x, y, t \in R$, and
\begin{equation}
\label{multiplication in T^+(3, R)}
\left(\begin{array}{ccc} e & x & t \\ 0 & e & y \\ 0 & 0 & e \end{array}\right)
\left(\begin{array}{ccc} e & x' & t' \\ 0 & e & y' \\ 0 & 0 & e 
                                                             \end{array}\right)
 = \left(\begin{array}{ccc} e & x + x' & t + t' + x \, y' \\
                               0 & e & y + y' \\ 0 & 0 & e \end{array}\right)
\end{equation}
for every $x, x', y, y', t, t' \in R$.  This implies that the mapping
from $(x, y, t) \in R^3$ to (\ref{H_1(R) into T^+(3, R)}) defines a
group isomorphism from the first Heisenberg group $H_1(R)$ onto
$T^+(3, R)$.  There is also a natural embedding of $H_n(R)$ into
$T^+(n + 2, R)$ for every positive integer $n$, which sends $(x, y, t)
\in H_n(R)$ to the following matrix $A = \{a_{j, k}\}_{j, k = 1}^{n + 2}$.
\begin{equation}
\label{H_n(R) into T^+(n + 2, R)}
\left(\begin{array}{ccccc} e      & x_1    & \cdots & x_n    & t \\
                           0      & e      & \cdots & 0      & y_1 \\
                           \vdots & \vdots & \ddots & \vdots & \vdots      \\
                           0      & 0      & \cdots & e      & y_n \\
                           0      & 0      & \cdots & 0      & e 
                                                    \end{array}\right)
\end{equation}
More precisely, $a_{j, k} = 0$ when $j > k$, $a_{j, j} = e$ for each
$j$, $a_{1, k + 1} = x_k$ for $k = 1, \ldots, n$, $a_{1, n + 2} = t$,
$a_{j + 1, n + 2} = y_j$ for $j = 1, \ldots, n$, and $a_{j, k} = 0$
when $2 \le j < k \le n + 1$.

        Let us restrict our attention to $n \ge 2$, since $T^+(1, R)$
is trivial.  If $B$ is any element of $M_n(R)$, then
\begin{equation}
\label{(I - B) (sum_{l = 0}^{n - 1} B^l) = ... = I - B^n}
        (I - B) \, \Big(\sum_{l = 0}^{n - 1} B^l\Big)
            = \Big(\sum_{l = 0}^{n - 1} B^l\Big) \, (I - B) = I - B^n,
\end{equation}
where $B^l = I$ when $l = 0$.  If $B \in T_n^+(R)$, then $B^n = 0$,
and it follows that $I - B$ is invertible in $M_n(R)$, with
\begin{equation}
\label{(I - B)^{-1} = sum_{l = 0}^{n - 1} B^l}
        (I - B)^{-1} = \sum_{l = 0}^{n - 1} B^l,
\end{equation}
which is contained in $T^+(n, R)$.  This gives a convenient expression
for the inverse of elements of $T^+(n, R)$, since every element of
$T^+(n, R)$ is of the form $I - B$ for some $B \in T_n^+(R)$.

        Let $r \in R$ be given, and let $\delta_r$ be the corresponding
ring homomorphism from $T_n(R)$ into itself discussed in the previous section.
Note that
\begin{equation}
\label{delta_r(A) in T^+(n, R)}
        \delta_r(A) \in T^+(n, R)
\end{equation}
for every $A \in T^+(n, R)$, so that $\delta_r$ defines a group
homomorphism from $T^+(n, R)$ into itself.  If $n = 2$, then
$\delta_r$ corresponds to the standard dilation on $R$ given by
multiplication by $r$, using the isomorphism between $R$ as a
commutative group with respect to addition and $T^+(2, R)$ indicated
by (\ref{R into T^+(2, R)}).  If $n = 3$, then $\delta_r$ corresponds
to (\ref{delta_r((x, y, t)) = (r x, r y, r^2 t)}) on $H_1(R)$, using
the isomorphism between $H_1(R)$ and $T^+(3, R)$ indicated by
(\ref{H_1(R) into T^+(3, R)}).  However, this does not work when $n
\ge 4$, using the embedding of $H_{n - 2}(R)$ into $T^+(n, R)$
indicated by (\ref{H_n(R) into T^+(n + 2, R)}).

        Now let $p$ be a prime number, so that the previous remarks
can be applied to $R = {\bf Q}_p$ and ${\bf Z}_p$.  Observe that
\begin{equation}
\label{T^+(n, {bf Z}_p) = GL(n, {bf Z}_p) cap T^+(n, {bf Q}_p)}
        T^+(n, {\bf Z}_p) = GL(n, {\bf Z}_p) \cap T^+(n, {\bf Q}_p)
\end{equation}
is a compact subgroup of $T^+(n, {\bf Q}_p)$ which is relatively open
in $T^+(n, {\bf Q}_p)$ as well.  Similarly, for each $l \in {\bf Z}_+$,
\begin{equation}
\label{GL_l(n, {bf Z}_p) cap T^+(n, {bf Z}_p)}
        GL_l(n, {\bf Z}_p) \cap T^+(n, {\bf Z}_p)
\end{equation}
is a compact relatively open subgroup of $T^+(n, {\bf Z}_p)$, and
hence of $T^+(n, {\bf Q}_p)$, consisting of $A = \{a_{j, k}\}_{j, k =
  1}^n \in T^+(n, {\bf Z}_p)$ such that $a_{j, k} \in p^l \, {\bf
  Z}_p$ when $j < k$.  Thus $T^+(n, {\bf Q}_p)$ has small compact open
subgroups, which are inherited from the corresponding subgroups of
$GL(n, {\bf Q}_p)$.  Because $GL_l(n, {\bf Z}_p)$ is a normal subgroup
of $GL(n, {\bf Z}_p)$ for each $l$, (\ref{GL_l(n, {bf Z}_p) cap T^+(n,
  {bf Z}_p)}) is a normal subgroup of $T^+(n, {\bf Z}_p)$ for each $l
\in {\bf Z}_+$, so that $T^+(n, {\bf Z}_p)$ has small compact open
normal subgroups.

        Let $N(A)$ be as in (\ref{N(A) = max_{1 le j < k le n} 
|a_{j, k}|_p^{1/(k - j)}}), where $A \in T_n({\bf Q}_p)$.  If we
restrict our attention to $A \in T^+(n, {\bf Q}_p)$, then $N(A) = 0$
if and only if $A = I$.  If $A, A' \in T^+(n, {\bf Q}_p)$, then
(\ref{N(A A') le max(N(A), N(A'))}) holds, since the diagonal entries
of $A$ and $A'$ are equal to $1$, which is an element of ${\bf Z}_p$.
Let us check that
\begin{equation}
\label{N(A^{-1}) le N(A)}
        N(A^{-1}) \le N(A)
\end{equation}
for every $A \in T^+(n, {\bf Q}_p)$, which implies that
\begin{equation}
\label{N(A^{-1}) = N(A)}
        N(A^{-1}) = N(A)
\end{equation}
for every $A \in T^+(n, {\bf Q}_p)$, by applying (\ref{N(A^{-1}) le
  N(A)}) to $A^{-1}$.  As before, we may as well take $A = I - B$,
where $B \in T_n^+({\bf Q}_p)$, so that $A^{-1}$ is given by (\ref{(I
  - B)^{-1} = sum_{l = 0}^{n - 1} B^l}).  Thus
\begin{equation}
\label{N(A^{-1}) = N((I - B)^{-1}) le max_{1 le l le n - 1} N(B^l)}
        N(A^{-1}) = N((I - B)^{-1}) \le \max_{1 \le l \le n - 1} N(B^l),
\end{equation}
by (\ref{N(A + A') le max(N(A), N(A'))}).  Clearly $N(A) = N(B)$, by
the definition (\ref{N(A) = max_{1 le j < k le n} |a_{j, k}|_p^{1/(k -
    j)}}) of $N(A)$, and
\begin{equation}
\label{N(B^l) le N(B)}
        N(B^l) \le N(B)
\end{equation}
for each $l \ge 1$, by repeated application of (\ref{N(A A') le
  max(N(A), N(A'))}).  This implies (\ref{N(A^{-1}) le N(A)}), as
desired.

        It follows that
\begin{equation}
\label{N((A')^{-1} A)}
        N((A')^{-1} \, A)
\end{equation}
defines a left-invariant ultrametric on $T^+(n, {\bf Q}_p)$, and that
\begin{equation}
\label{N(A (A')^{-1})}
        N(A \, (A')^{-1})
\end{equation}
defines a right-invariant ultrametric on $T^+(n, {\bf Q}_p)$, as in
Section \ref{topological groups}.  It is easy to see that the
topologies on $T^+(n, {\bf Q}_p)$ determined by these ultrametrics are
the same as the one induced by the usual topology on $M_n({\bf Q}_p)$.
These ultrametrics also behave well with respect to the dilations
$\delta_r$ on $T^+(n, {\bf Q}_p)$, by (\ref{N(delta_r(A)) = |r|_p
  N(A)}).  

        Note that Haar measure on $T^+(n, {\bf Q}_p)$ corresponds to
the product of Haar measure on $n (n - 1) / 2$ copies of ${\bf Q}_p$.
If $E \subseteq T^+(n, {\bf Q}_p)$ is a Borel set and $r \in {\bf
  Q}_p$, then $\delta_r(E)$ is a Borel set in $T^+(n, {\bf Q}_p)$ too, 
and the Haar measure of $\delta_r(E)$ is equal to
\begin{equation}
\label{|r|_p^{d(n)}}
        |r|_p^{d(n)}
\end{equation}
times the Haar measure of $E$, where
\begin{equation}
\label{d(n) = sum_{1 le j < k le n} (k - j) = sum_{l = 1}^{n - 1} (n - l) l}
        d(n) = \sum_{1 \le j < k \le n} (k - j) = \sum_{l = 1}^{n - 1} (n - l) \, l.
\end{equation}
This uses the analogous property of Haar measure on ${\bf Q}_p$, and
it implies that Haar measure on $T^+(n, {\bf Q}_p)$ is Ahlfors regular
of dimension $d(n)$ with respect to (\ref{N((A')^{-1} A)}) or
(\ref{N(A (A')^{-1})}).

        If $l$ is any integer, then
\begin{eqnarray}
\label{{A in T^+(n, {bf Q}_p) : N(A) le p^{-l}} = ...}
\lefteqn{\{A \in T^+(n, {\bf Q}_p) : N(A) \le p^{-l}\}} \\
 & = & \{A \in T^+(n, {\bf Q}_p) : a_{j, k} \in p^{l \, (k - j)} \, {\bf Z}_p
                                    \hbox{ when } j < k\} \nonumber\\
 & = & \delta_{p^l}(T^+(n, {\bf Z}_p))                      \nonumber
\end{eqnarray}
is a compact open subgroup of $T^+(n, {\bf Q}_p)$, so that $T^+(n,
{\bf Z}_p)$ has large compact open subgroups in particular.
Similarly, if $l_1, \ldots, l_{n - 1}$ are integers such that
\begin{equation}
\label{l_{alpha + beta} le l_alpha + l_beta}
        l_{\alpha + \beta} \le l_\alpha + l_\beta
\end{equation}
for every $\alpha, \beta \ge 1$ with $\alpha + \beta \le n - 1$, then
one can check that
\begin{equation}
\label{{A in T^+(n, {bf Q}_p) : a_{j, k} in p^{l_{k - j}} {bf Z}_p when j < k}}
        \{A \in T^+(n, {\bf Q}_p) : a_{j, k} \in p^{l_{k - j}} \, {\bf Z}_p
                                             \hbox{ when } j < k\}
\end{equation}
is a compact open subgroup of $T^+(n, {\bf Q}_p)$.  More precisely,
\begin{equation}
\label{{B in T_n^+({bf Q}_p) : b_{j, k} in p^{l_{k - j}} {bf Z}_p when j < k}}
        \{B \in T_n^+({\bf Q}_p) : b_{j, k} \in p^{l_{k - j}} \, {\bf Z}_p
                                            \hbox{ when } j < k\}
\end{equation}
is a compact open sub-ring of $T_n^+({\bf Q}_p)$, and (\ref{{A in
    T^+(n, {bf Q}_p) : a_{j, k} in p^{l_{k - j}} {bf Z}_p when j <
    k}}) is the same as the collection of $A \in T^+(n, {\bf Q}_p)$
such that $A - I$ is in (\ref{{B in T_n^+({bf Q}_p) : b_{j, k} in
    p^{l_{k - j}} {bf Z}_p when j < k}}).

\section{Upper-triangular matrices, 3}
\label{upper-triangular matrices, 3}
\setcounter{equation}{0}

        Let $R$ be a commutative ring with nonzero multiplicative
identity element $e$, and let $R^*$ be the multiplicative group of
invertible elements in $R$.  Also let $n$ be a positive integer, and
let $M_n(R)$, $T_n(R)$, $T_n^+(R)$, and $T^+(n, R)$ be as in the
previous sections.  Consider the collection $T(n, R)$ of $A = \{a_{j,
  k}\}_{j, k = 1}^n$ in $T_n(R)$ whose diagonal entries are invertible
in $R$, so that $a_{j, j} \in R^*$ for each $j = 1, \ldots, n$.  Thus
\begin{equation}
\label{det A = prod_{j = 1}^n a_{j, j}}
        \det A = \prod_{j = 1}^n a_{j, j}
\end{equation}
is invertible in $R$ when $A \in T(n, R)$.  Conversely, if $A \in
T_n(R)$ and $\det A$ is invertible in $R$, then $a_{j, j} \in R^*$ for
each $j = 1, \ldots, n$, so that $A \in T(n, R)$.

        Equivalently,
\begin{equation}
\label{T(n, R) = T_n(R) cap GL(n, R)}
        T(n, R) = T_n(R) \cap GL(n, R).
\end{equation}
In particular, if $A, A' \in T(n, R)$, then their product $A \, A'$ is
in $T(n, R)$, because of the analogous statements for $T_n(R)$ and
$GL(n, R)$.  Alternatively, the diagonal entries of $A \, A'$ are
equal to the products of the corresponding diagonal entries of $A$ and
$A'$, and hence are invertible in $R$ when the diagonal entries of $A$
and $A'$ are invertible in $R$.  If $A \in T(n, R)$, then it is well
known that $A^{-1} \in T(n, R)$, so that $T(n, R)$ is a subgroup of
$GL(n, R)$.  One way to see this is to express $A$ as the product of a
diagonal matrix with the same diagonal entries as $A$ and an element
of $T^+(n, R)$, and then invert the two factors.

        Note that the collection of diagonal $n \times n$ matrices with
entries in $R$ is a sub-ring of $T_n(R)$ and hence $M_n(R)$, and that
the mapping from $A \in T_n(R)$ to the diagonal matrix with the same
diagonal entries as $A$ is a ring homomorphism.  The restriction of
this mapping to $T(n, R)$ defines a group homomorphism from $T(n, R)$
onto the group of diagonal matrices with diagonal entries in $R^*$,
and $T^+(n, R)$ is the kernel of this homomorphism. In particular,
$T^+(n, R)$ is a normal subgroup of $T(n, R)$.

        Now let $p$ be a prime number, so that the previous remarks
can be applied to $R = {\bf Q}_p$ and ${\bf Z}_p$.  Thus $T(n, {\bf
  Q}_p)$ is a relatively open set in $T_n({\bf Q}_p)$ with respect to
the topology induced by the usual one on $M_n({\bf Q}_p)$, and a
relatively closed subgroup of $GL(n, {\bf Q}_p)$.  By definition,
\begin{equation}
\label{T(n, Z_p) = {A in T_n(Z_p) : |a_{j, j}|_p = 1 for each j = 1, ldots, n}}
        T(n, {\bf Z}_p) = \{A \in T_n({\bf Z}_p) : |a_{j, j}|_p = 1 
                                 \hbox{ for each } j = 1, \ldots, n\},
\end{equation}
since $x \in {\bf Z}_p$ is invertible in ${\bf Z}_p$ if and only if
$|x|_p = 1$.  Equivalently,
\begin{equation}
\label{T(n, {bf Z}_p) = {A in T_n({bf Z}_p) : |det A|_p = 1}}
        T(n, {\bf Z}_p) = \{A \in T_n({\bf Z}_p) : |\det A|_p = 1\}.
\end{equation}
More precisely, if $A \in T(n, {\bf Z}_p)$, then $|a_{j, j}|_p = 1$
for $j = 1, \ldots, n$, and hence $|\det A|_p = 1$, by (\ref{det A =
  prod_{j = 1}^n a_{j, j}}).  Similarly, if $A \in T_n(A)$, then
$a_{j, j} \in {\bf Z}_p$ for $j = 1, \ldots, n$, and $\det A \in {\bf Z}_p$.
In this case, if we also have that $|\det A|_p = 1$, then it follows
that $|a_{j, j}|_p = 1$ for each $j$, so that $A \in T(n, {\bf Z}_p)$.

        Of course, $T(n, {\bf Z}_p)$ is a compact subgroup of
$GL(n, {\bf Z}_p)$, and $T(n, {\bf Z}_p)$ is relatively open in 
$T(n, {\bf Q}_p)$.  Similarly,
\begin{equation}
\label{GL_l(n, {bf Z}_p) cap T(n, {bf Z}_p)}
        GL_l(n, {\bf Z}_p) \cap T(n, {\bf Z}_p)
\end{equation}
is a compact relatively open normal subgroup of $T(n, {\bf Z}_p)$ for
each positive integer $l$, and hence a compact relatively open
subgroup of $T(n, {\bf Q}_p)$.  This shows that $T(n, {\bf Z}_p)$ has
small compact open normal subgroups, and that $T(n, {\bf Q}_p)$ has
small compact open subgroups.  However, $T(n, {\bf Q}_p)$ does not
have large compact subgroups, for essentially the same reasons as
for ${\bf Q}_p^*$.

        As a variant of this, consider
\begin{equation}
\label{widetilde{T}(n, Q_p) = {A in T(n, Q_p) : |a_{j, j}|_p = 1}}
        \widetilde{T}(n, {\bf Q}_p) = \{A \in T(n, {\bf Q}_p) : |a_{j, j}|_p = 1
                                   \hbox{ for each } j = 1, \ldots, n\},
\end{equation}
which is a relatively open normal subgroup of $T(n, {\bf Q}_p)$ that
contains $T(n, {\bf Z}_p)$ as a relatively open compact subgroup.
Thus $\widetilde{T}(n, {\bf Q}_p)$ also contains (\ref{GL_l(n, {bf
    Z}_p) cap T(n, {bf Z}_p)}) as a compact relatively open subgroup
for each $l \in {\bf Z}_+$, so that $\widetilde{T}(n, {\bf Q}_p)$ has
small compact open subgroups as well.  By construction, $T^+(n, {\bf
  Q}_p)$ is a relatively closed normal subgroup of $\widetilde{T}(n,
{\bf Q}_p)$, and every element of $\widetilde{T}(n, {\bf Q}_p)$ can be
expressed as the product of a diagonal matrix whose diagonal entries
have $p$-adic absolute value equal to $1$ and an element of $T^+(n,
{\bf Q}_p)$.

        Let $N(A)$ be as in (\ref{N(A) = max_{1 le j < k le n} 
|a_{j, k}|_p^{1/(k - j)}}) for each $A \in T_n({\bf Q}_p)$, and for
$A \in \widetilde{T}(n, {\bf Q}_p)$ in particular.  Thus $N(A) = 0$ exactly
when $A$ is a diagonal matrix.  Note that (\ref{N(A A') le max(N(A),
  N(A'))}) holds when $A, A' \in \widetilde{T}(n, {\bf Q}_p)$, since
the diagonal entries of $A$ and $A'$ are in ${\bf Z}_p$.  We have
already seen that (\ref{N(A^{-1}) le N(A)}) and hence (\ref{N(A^{-1})
  = N(A)}) hold for every $A \in T^+(n, {\bf Q}_p)$, and this can be
extended to $A \in \widetilde{T}(n, {\bf Q}_p)$ by expressing $A$ as
the product of a diagonal matrix and an element of $T^+(n, {\bf
  Q}_p)$.  It follows that (\ref{N((A')^{-1} A)}) defines a
left-invariant semi-ultrametric on $\widetilde{T}(n, {\bf Q}_p)$, and
that (\ref{N(A (A')^{-1})}) defines a right-invariant semi-ultrametric
on $\widetilde{T}(n, {\bf Q}_p)$.  It is easy to see that these
semi-ultrametrics are compatible with the topology on
$\widetilde{T}(n, {\bf Q}_p)$ induced by the standard topology on
$M_n({\bf Q}_p)$, and that they are proper.  This implies that
$\widetilde{T}(n, {\bf Q}_p)$ has large compact open subgroups, as
before.

        Of course,
\begin{equation}
\label{D(A, A') = max_{1 le j le n} |a_{j, j} - a'_{j, j}|_p}
        D(A, A') = \max_{1 \le j \le n} |a_{j, j} - a'_{j, j}|_p
\end{equation}
is also a semi-ultrametric on $\widetilde{T}(n, {\bf Q}_p)$, which is
equal to $0$ exactly when $A$ and $A'$ have the same diagonal entries.
This semi-ultrametric is invariant under both left and right translations
on $\widetilde{T}(n, {\bf Q}_p)$ as a group with respect to matrix
multiplication, because $A \mapsto a_{j, j}$ is a homomorphism from
$\widetilde{T}(n, {\bf Q}_p)$ onto the multiplicative group of $p$-adic 
numbers with $p$-adic absolute value equal to $1$ for each $j$.  
The maximum of (\ref{N((A')^{-1} A)}) and (\ref{D(A, A') = max_{1 le j
    le n} |a_{j, j} - a'_{j, j}|_p}) is a left-invariant ultrametric
on $\widetilde{T}(n, {\bf Q}_p)$ that determines the usual topology on
$\widetilde{T}(n, {\bf Q}_p)$, and the maximum of (\ref{N(A
  (A')^{-1})}) and (\ref{D(A, A') = max_{1 le j le n} |a_{j, j} -
  a'_{j, j}|_p}) is a right-invariant ultrametric on $\widetilde{T}(n,
{\bf Q}_p)$ that determines the same topology on $\widetilde{T}(n,
{\bf Q}_p)$.

\section{$a \, x + b$ Groups}
\label{a x + b groups}
\setcounter{equation}{0}

        Let $R$ be a commutative ring with a nonzero multiplicative identity 
element $e$ again, and let $A(R)$ be the collection of ``affine'' mappings
on $R$, which is to say mappings from $R$ into itself of the form
\begin{equation}
\label{f(x) = a x + b}
        f(x) = a \, x + b
\end{equation}
for some $a, b \in R$.  If
\begin{equation}
\label{g(x) = c x + d}
        g(x) = c \, x + d
\end{equation}
is another element of $A(R)$ for some $c, d \in R$, then the composition
\begin{equation}
\label{(f circ g)(x) = f(g(x)) = a c x + a d + b}
        (f \circ g)(x) = f(g(x)) = a \, c \, x + a \, d + b
\end{equation}
of $f$ and $g$ is also in $A(R)$.  It follows that $A(R)$ is a
semigroup with respect to composition, and the identity mapping on $R$
is the identity element in $A(R)$.  If $f \in A(R)$ is of the form
(\ref{f(x) = a x + b}), then $a, b \in R$ are uniquely determined by
$f$ as a mapping from $R$ into itself, because
\begin{equation}
\label{a = f(e) - f(0) and b = f(0)}
        a = f(e) - f(0) \quad\hbox{and}\quad b = f(0).
\end{equation}
Thus we can identify $A(R)$ as a set with $R \times R$.  

        Let $R^*$ be the multiplicative group of invertible elements
of $R$, and let $A^*(R)$ be the collection of $f \in A(R)$ of the form
(\ref{f(x) = a x + b}) with $a \in R^*$.  If $f \in A^*(R)$, then it is easy
to see that $f$ is an invertible mapping from $R$ onto itself, with
\begin{equation}
\label{f^{-1}(x) = a^{-1} x - a^{-1} b}
        f^{-1}(x) = a^{-1} \, x - a^{-1} \, b,
\end{equation}
where $a, b \in R$ are as in (\ref{f(x) = a x + b}) again.
Conversely, if $f \in A(R)$ maps $R$ onto itself, then there is an $x
\in R$ such that $f(x) = e + b$, which implies that $a \in R^*$.  Thus
$A^*(R)$ is the same as the collection of $f \in A(R)$ such that $f$
is an invertible mapping from $R$ onto itself, in which case $f^{-1}
\in A^*(R)$ too.  The composition of two elements of $A^*(R)$ is also
an element of $A^*(R)$, so that $A^*(R)$ is a group with respect to
composition of mappings.  

        Let $U$ be a subset of $R$, and let $A(U, R)$ be the collection 
of $f \in A(R)$ of the form (\ref{f(x) = a x + b}) with $a \in U$.  If
$U$ is a sub-semigroup of $R$ with respect to multiplication, then
$A(U, R)$ is a sub-semigroup of $A(R)$ with respect to composition of
mappings.  Similarly, if $U$ is a subgroup of $R^*$, then $A(U, R)$ is
a subgroup of $A^*(R) = A(R^*, R)$.  In both cases, there is a natural
homomorphism from $A(U, R)$ onto $U$, which sends $f \in A(U, R)$ as
in (\ref{f(x) = a x + b}) to $a \in U$.  If $e \in U$, which is
automatic when $U$ is a subgroup of $R^*$, then the kernel of this
homomorphism is equal to $A(\{e\}, R)$, which consists of the
translation mappings $f(x) = x + b$ with $b \in R$.

        If $R_1$ is a subring of $R$ that contains $e$, then we can identify 
$A(R_1)$ with a sub-semigroup of $A(R)$, consisting of the $f \in A(R)$ 
as in (\ref{f(x) = a x + b}) with $a, b \in R_1$.  Equivalently, this
happens exactly when $f(R_1) \subseteq R_1$, because of (\ref{a = f(e)
  - f(0) and b = f(0)}).  In this case, $R_1^*$ is a subgroup of
$R^*$, and $A^*(R_1)$ can be identified with a subgroup of $A^*(R)$.
However, an element of $R_1$ may be invertible as an element of $R$,
without its inverse being in $R_1$.  This means that an element of
$A(R_1)$ may correspond to an element of $A^*(R)$, without being in
$A^*(R_1)$.

        If $R$ is a topological ring, then we can get a topology on
$A(R)$ by identifying $A(R)$ with $R \times R$ as before, and using
the product topology on $R \times R$.  It is easy to see that the
mapping from $f, g \in A(R)$ to $f \circ g$ is continuous as a mapping
from $A(R) \times A(R)$ into $A(R)$ with respect to the product topology
on $A(R)$, so that $A(R)$ becomes a topological semigroup under these 
conditions.  If $R$ is a topological field, then $R^* = R \setminus \{0\}$
is an open set in $R$, and hence $A^*(R)$ is an open set in $A(R)$.
In this case, the mapping from $f \in A^*(R)$ to $f^{-1}$ is also continuous,
so that $A^*(R)$ is a topological group.

        Note that the mapping 
\begin{equation}
\label{b in R mapsto f(x) = x + b in A(R)}
        b \in R \mapsto f(x) = x + b \in A(R)
\end{equation}
defines an isomorphism from $R$ as a commutative group with respect to
addition onto $A(\{e\}, R)$.  Similarly,
\begin{equation}
\label{a in R mapsto f(x) = a x in A(R)}
        a \in R \mapsto f(x) = a \, x \in A(R)
\end{equation}
defines an isomorphism from $R$ as a semigroup with respect to
multiplication onto a sub-semigroup of $A(R)$.  The restriction of
(\ref{a in R mapsto f(x) = a x in A(R)}) to $a \in R^*$ defines an
isomorphism from $R^*$ onto a subgroup of $A^*(R)$.  If $R$ is a
topological ring, then these mappings are homeomorphisms onto their
images in $A(R)$, with respect to the induced topology.

        We can also think of $A(R)$ as a module over $R$, 
with respect to pointwise addition and scalar multiplication of affine
functions on $R$.  If we identify $A(R)$ as a set with $R \times R$ as
before, then this corresponds exactly to coordinatewise addition and
scalar multiplication on $R \times R$.  In particular, if $R$ is a
field, then $A(R)$ may be considered as a $2$-dimensional vector space
over $R$.  If $R$ is a topological ring or field, then $A(R)$ is a topological
module or vector space over $R$, as appropriate, in the sense that addition
and scalar multiplication are continuous.

        There is an obvious identification between $f \in A(R)$ as in 
(\ref{f(x) = a x + b}) and the $2 \times 2$ matrix
\begin{equation}
\label{( a b / 0 e )}
\left(\begin{array}{cc} a & b \\ 0 & e \end{array}\right).
\end{equation}
In particular, the identity mapping $f(x) = x$ corresponds to the
identity matrix $I = \left({e \atop 0}{0 \atop e}\right)$, and
\begin{equation}
\label{( a b / 0 e) (c d / 0 e) = ...}
\left(\begin{array}{cc} a & b \\ 0 & e \end{array}\right)
 \left(\begin{array}{cc} c & d \\ 0 & e \end{array}\right)
  = \left(\begin{array}{cc} a \, c & a \, d + b \\ 0 & e \end{array}\right)
\end{equation}
for every $a, b, c, d \in R$, as in (\ref{(f circ g)(x) = f(g(x)) = a
  c x + a d + b}).  If $a \in R^*$ and $b \in R$, then (\ref{( a b / 0
  e )}) is invertible as a $2 \times 2$ matrix with entries in $R$,
and its inverse is the matrix of the same form corresponding to
(\ref{f^{-1}(x) = a^{-1} x - a^{-1} b}).  Thus $A(R)$ corresponds to a
sub-semigroup of $T_2(R)$ as in Section \ref{upper-triangular
  matrices} with respect to matrix multiplication, $A(\{e\}, R)$
corresponds to $T^+(2, R)$ as in Section \ref{upper-triangular
  matrices, 2}, and $A^*(R)$ corresponds to a subgroup of $T(2, R)$ as
in Section \ref{upper-triangular matrices, 3}, with respect to this
identification.

        Now let $p$ be a prime number, so that the previous remarks 
can be applied to $R = {\bf Q}_p$.  In particular, $A({\bf Q}_p)$ may
be considered as a $2$-dimensional vector space over ${\bf Q}_p$, with
respect to pointwise addition and scalar multiplication of affine
functions on ${\bf Q}_p$.  If $f \in A({\bf Q}_p)$ is as in (\ref{f(x)
  = a x + b}), then put
\begin{equation}
\label{||f||_{A({bf Q}_p)} = max(|a|_p, |b|_p)}
        \|f\|_{A({\bf Q}_p)} = \max(|a|_p, |b|_p),
\end{equation}
which defines an ultranorm on $A({\bf Q}_p)$ as a vector space over
${\bf Q}_p$.  If we identify $A({\bf Q}_p)$ with ${\bf Q}_p \times
{\bf Q}_p$ as a two-dimensional vector space over ${\bf Q}_p$ in the
usual way, then $\|f\|_{A({\bf Q}_p)}$ is the same as the standard
ultranorm on ${\bf Q}_p \times {\bf Q}_p \cong {\bf Q}_p^2$.  Of
course, the topology on $A({\bf Q}_p)$ determined by the corresponding
ultrametric
\begin{equation}
\label{||f - g||_{A({bf Q}_p)}}
        \|f - g\|_{A({\bf Q}_p)}
\end{equation}
is the same as the one that we get from the product topology on ${\bf
  Q}_p \times {\bf Q}_p$.

        The space $M_2({\bf Q}_p)$ of $2 \times 2$ matrices with entries 
in ${\bf Q}_p$ is a four-dimensional vector space over ${\bf Q}_p$ in
the usual way, and the matrices of the form (\ref{( a b / 0 e )}) span
a two-dimensional affine subspace of $M_2({\bf Q}_p)$.  As in Section
\ref{ultranorms}, the standard ultranorm of an element of $M_2({\bf
  Q}_p)$ is defined by taking the maximum of the $p$-adic absolute
values of its entries.  If $a, b \in {\bf Q}_p$ and $f \in A({\bf
  Q}_p)$ is as in (\ref{f(x) = a x + b}), then the norm of the
corresponding matrix (\ref{( a b / 0 e )}) is equal to
\begin{equation}
\label{max(|a|_p, |b|_p, 1) = max(||f||_{A({bf Q}_p)}, 1)}
        \max(|a|_p, |b|_p, 1) = \max(\|f\|_{A({\bf Q}_p)}, 1).
\end{equation}
However, if $g \in A({\bf Q}_p)$ is as in (\ref{g(x) = c x + d}), then
the difference between the matrices corresponding to $f$ and $g$ is
equal to
\begin{equation}
\label{{a b / 0 e} - {c d / 0 e} = ...}
\left(\begin{array}{cc} a & b \\ 0 & e \end{array}\right)
 - \left(\begin{array}{cc} c & d \\ 0 & e \end{array}\right)
 = \left(\begin{array}{cc} a - c & b - d \\ 0 & 0 \end{array}\right),
\end{equation}
and the norm of (\ref{{a b / 0 e} - {c d / 0 e} = ...}) is equal to
(\ref{||f - g||_{A({bf Q}_p)}}).

        Let $U_p$ be the multiplicative group of $x \in {\bf Q}_p$
with $|x|_p = 1$, which is the same as the group of invertible
elements in ${\bf Z}_p$.  This is a compact open subgroup of ${\bf
  Q}_p^* = {\bf Q}_p \setminus \{0\}$, and $A(U_p, {\bf Q}_p)$ is an
open normal subgroup of $A^*({\bf Q}_p)$.  Note that the elements of
$A(U_p, {\bf Q}_p)$ are isometries on ${\bf Q}_p$ with respect to the
$p$-adic metric.  If $f \in A({\bf Q}_p)$, then
\begin{equation}
\label{||f||_{A({bf Q}_p)} = max {|f(x)|_p : x in {bf Z}_p}}
        \|f\|_{A({\bf Q}_p)} = \max \{|f(x)|_p : x \in {\bf Z}_p\},
\end{equation}
because $|f(x)|_p \le \|f\|_{A({\bf Q}_p)}$ for every $x \in {\bf
  Z}_p$, and
\begin{equation}
\label{||f||_{A({bf Q}_p)} le max(|f(0)|_p, |f(1)|_p)}
        \|f\|_{A({\bf Q}_p)} \le \max(|f(0)|_p, |f(1)|_p),
\end{equation}
by (\ref{a = f(e) - f(0) and b = f(0)}).  If $\alpha \in A(U_p, {\bf
  Q}_p)$, then it follows that
\begin{equation}
\label{||alpha circ f - alpha circ g||_{A({bf Q}_p)} = ||f - g||_{A({bf Q}_p)}}
 \|\alpha \circ f - \alpha \circ g\|_{A({\bf Q}_p)} = \|f - g\|_{A({\bf Q}_p)}
\end{equation}
for every $f, g \in A({\bf Q}_p)$.  Of course, this can also be
verified using (\ref{||f||_{A({bf Q}_p)} = max(|a|_p, |b|_p)}).  In
particular, the restriction of (\ref{||f - g||_{A({bf Q}_p)}}) to $f$,
$g$ in $A(U_p, {\bf Q}_p)$ is invariant under left translations on
$A(U_p, {\bf Q}_p)$, as a group with respect to composition of
mappings.  If $f \in A({\bf Q}_p)$ as in (\ref{f(x) = a x + b}) is
identified with the $2 \times 2$ matrix (\ref{( a b / 0 e )}), then
$A(U_p, {\bf Q}_p)$ corresponds to a closed subgroup of
$\widetilde{T}(2, {\bf Q}_p)$ as in Section \ref{upper-triangular
  matrices, 3}.

        As before, we can identify $A({\bf Z}_p)$ with the set of 
$f \in A({\bf Q}_p)$ with coefficients in ${\bf Z}_p$, which is a compact
open sub-semigroup of $A({\bf Q}_p)$.  Similarly, $A^*({\bf Z}_p) =
A(U_p, {\bf Z}_p)$ corresponds to a compact open subgroup of $A^*({\bf
  Q}_p)$.  Using (\ref{||f||_{A({bf Q}_p)} = max {|f(x)|_p : x in {bf Z}_p}}), 
we get that
\begin{equation}
\label{||f circ beta - g circ beta||_{A({bf Q}_p)} = ||f - g||_{A({bf Q}_p)}}
        \|f \circ \beta - g \circ \beta\|_{A({\bf Q}_p)} = \|f - g\|_{A({\bf Q}_p)}
\end{equation}
for every $\beta \in A^*({\bf Z}_p)$ and $f, g \in A({\bf Q}_p)$,
because $\beta({\bf Z}_p) = {\bf Z}_p$.  This implies that the
restriction of (\ref{||f - g||_{A({bf Q}_p)}}) to $A^*({\bf Z}_p)$ is
invariant under both left and right translations.  Of course, one can
also look at this in terms of matrices, since $A^*({\bf Z}_p)$
corresponds to a subgroup of $T(2, {\bf Z}_p)$, which is contained in
$GL(2, {\bf Z}_p)$.

        If $f \in A(U_p, {\bf Q}_p)$ is as in (\ref{f(x) = a x + b}), then
one can check that
\begin{equation}
\label{L(f) = ||f(x) - x||_{A({bf Q}_p)} = max(|a - 1|_p, |b|_p)}
        L(f) = \|f(x) - x\|_{A({\bf Q}_p)} = \max(|a - 1|_p, |b|_p)
\end{equation}
satisfies
\begin{equation}
\label{L(f^{-1}) = L(f)}
        L(f^{-1}) = L(f),
\end{equation}
using (\ref{f^{-1}(x) = a^{-1} x - a^{-1} b}) or (\ref{||alpha circ f
  - alpha circ g||_{A({bf Q}_p)} = ||f - g||_{A({bf Q}_p)}}).
Similarly, if $f, g \in A(U_p, {\bf Q}_p)$, then
\begin{eqnarray}
\label{L(f circ g) = ...  = max(L(f), L(g))}
        L(f \circ g) & = & \|f(g(x)) - x\|_{A({\bf Q}_p)}          \\
      & \le & \max(\|f(g(x)) - f(x)\|_{A({\bf Q}_p)}, \|f(x) - x\|_{A({\bf Q}_p)})
                                                       \nonumber \\
 & = & \max(\|g(x) - x\|_{A({\bf Q}_p)}, \|f(x) - x\|_{A({\bf Q}_p)}) \nonumber \\
 & = & \max(L(f), L(g)), \nonumber
\end{eqnarray}
by (\ref{||alpha circ f - alpha circ g||_{A({bf Q}_p)} = ||f -
  g||_{A({bf Q}_p)}}).  Let $t$ be a real number with $t \ge 1$, and put
\begin{equation}
\label{L'(f) = L(f) le 1}
        L'(f) = L(f) \le 1
\end{equation}
when $f \in A^*({\bf Z}_p)$, and
\begin{equation}
\label{L'(f) = t}
        L'(f) = t
\end{equation}
when $f \in A^*({\bf Q}_p)$ is not in $A^*({\bf Z}_p)$.  Thus
\begin{equation}
\label{L'(f^{-1}) = L'(f)}
        L'(f^{-1}) = L'(f)
\end{equation}
for every $f \in A^*({\bf Q}_p)$, and
\begin{equation}
\label{L'(f circ g) le max(L'(f), L'(g))}
        L'(f \circ g) \le \max(L'(f), L'(g))
\end{equation}
for every $f, g \in A^*({\bf Q}_p)$, by the previous statements for $L$.
It follows that
\begin{equation}
\label{L'(g^{-1} circ f)}
        L'(g^{-1} \circ f)
\end{equation}
defines an ultrametric on $A^*({\bf Q}_p)$ that is invariant under
left translations on $A^*({\bf Q}_p)$, as a group with respect to
composition of mappings.  If $f, g \in A^*({\bf Z}_p)$, then
(\ref{L'(g^{-1} circ f)}) is equal to $L(g^{-1} \circ f)$, which is
the same as (\ref{||f - g||_{A({bf Q}_p)}}), by (\ref{||alpha circ f -
  alpha circ g||_{A({bf Q}_p)} = ||f - g||_{A({bf Q}_p)}}).  This
implies that the topology on $A^*({\bf Q}_p)$ determined by
(\ref{L'(g^{-1} circ f)}) is the same as the usual topology on
$A^*({\bf Q}_p)$ discussed earlier, since $A^*({\bf Z}_p)$ is an open
subgroup of $A^*({\bf Q}_p)$.  Note that $A^*({\bf Q}_p)$ has small
compact open subgroups, and that $A^*({\bf Z}_p)$ has small open
normal subgroups.  Similarly, $A(U_p, {\bf Q}_p)$ has large compact
open subgroups, and the restriction of the ultrametric (\ref{||f -
  g||_{A({bf Q}_p)}}) to $A(U_p,{\bf Q}_p)$ is proper.  However,
$A^*({\bf Q}_p)$ does not have large compact subgroups, basically
because ${\bf Q}_p^*$ does not have large compact subgroups.

\section{Cells in ${\bf Q}_p$}
\label{cells in Q_p}
\setcounter{equation}{0}

        Let $p$ be a prime number, and let $A^*({\bf Q}_p)$ be as in
the previous section.  By a \emph{cell} in ${\bf Q}_p$ we mean a
subset of ${\bf Q}_p$ of the form $f({\bf Z}_p)$ for some $f \in
A^*({\bf Q}_p)$.  Equivalently, $C \subseteq {\bf Q}_p$ is a cell if
$C$ is a closed ball in ${\bf Q}_p$ of positive radius with respect to
the $p$-adic metric.  Let $\mathcal{C}({\bf Q}_p)$ be the collection
of all cells in ${\bf Q}_p$.

        Let us say that a cell $C \in \mathcal{C}({\bf Q}_p)$ is the
child of a cell $C' \in \mathcal{C}({\bf Q}_p)$ if $C \subseteq C'$
and the diameter of $C$ is equal to $1/p$ times the diameter of $C'$.
Each cell in ${\bf Q}_p$ has exactly $p$ children in $\mathcal{C}({\bf Q}_p)$, 
and is the child of exactly one cell in $\mathcal{C}({\bf Q}_p)$.
This leads to a natural graph whose vertices are cells in ${\bf Q}_p$,
and for which two cells are adjacent when one is a child of the other.
It is easy to see that this graph is connected, because any two cells in 
${\bf Q}_p$ are contained in another cell.  In fact, this graph is a tree.

        If $f \in A^*({\bf Q}_p)$ and $C$ is a cell in ${\bf Q}_p$,
then $f(C)$ is a cell in ${\bf Q}_p$ too, and
\begin{equation}
\label{C mapsto f(C)}
        C \mapsto f(C)
\end{equation}
defines a one-to-one mapping from $\mathcal{C}({\bf Q}_p)$ onto
itself.  This mapping actually leads to an automorphism on the tree
associated to $\mathcal{C}({\bf Q}_p)$, because it sends the children
of a cell $C$ to the children of $f(C)$.  Thus we get a homomorphism
from $A^*({\bf Q}_p)$ into the automorphism group of this tree.  If
$f$ and $g$ are distinct elements of $A^*({\bf Q}_p)$, then one can
check that there is a cell $C$ in ${\bf Q}_p$ such that $f(C) \ne
g(C)$, so that this homomorphism is an embedding.

        Let $\rho(C, C')$ be the path metric on $\mathcal{C}({\bf Q}_p)$,
which is the length of the shortest path between two cells $C$, $C'$
in the corresponding graph, where the length of each edge is equal to $1$.
If $C \subseteq C'$, then this path consists of the sequence of
cells containing $C$ and contained in $C'$.  Otherwise, if $C''$ is the
smallest cell that contains both $C$ and $C'$, then the shortest
path from $C$ to $C'$ can be obtained by combining the minimal paths between
$C$ and $C'$ to $C''$.  Note that
\begin{equation}
\label{rho(f(C), f(C')) = rho(C, C')}
        \rho(f(C), f(C')) = \rho(C, C')
\end{equation}
for every $f \in A^*({\bf Q}_p)$ and $C, C' \in \mathcal{C}({\bf
  Q}_p)$, since the mapping (\ref{C mapsto f(C)}) determines an
automorphism of the graph corresponding to $\mathcal{C}({\bf Q}_p)$.

        As before, $A^*({\bf Z}_p)$ may be identified with a compact open
subgroup of $A^*({\bf Q}_p)$, consisting of the $f \in A^*({\bf Q}_p)$
such that $f({\bf Z}_p) = {\bf Z}_p$.  Although $A^*({\bf Z}_p)$ is
not a normal subgroup of $A^*({\bf Q}_p)$, one can still consider the
quotient space $A^*({\bf Q}_p) / A^*({\bf Z}_p)$ as a set, with the
natural action of $A^*({\bf Q}_p)$ on the left.  The mapping from $f
\in A^*({\bf Q}_p)$ to $f({\bf Z}_p) \in \mathcal{C}({\bf Q}_p)$ leads
to a one-to-one mapping from $A^*({\bf Q}_p) / A^*({\bf Z}_p)$ onto
$\mathcal{C}({\bf Q}_p)$.  By construction, the natural action of
$A^*({\bf Q}_p)$ on $A^*({\bf Q}_p) / A^*({\bf Z}_p)$ on the left
corresponds exactly to the action of $A^*({\bf Q}_p)$ on
$\mathcal{C}({\bf Q}_p)$ already defined.

        Observe that
\begin{equation}
\label{rho(f({bf Z}_p), g({bf Z}_p))}
        \rho(f({\bf Z}_p), g({\bf Z}_p))
\end{equation}
defines a semimetric on $A^*({\bf Q}_p)$ that is invariant under left
translations on $A^*({\bf Q}_p)$, as a group with respect to
composition of mappings.  This semimetric is equal to $0$ when $f({\bf
  Z}_p) = g({\bf Z}_p)$, which is to say that $g^{-1} \circ f \in
A^*({\bf Z}_p)$.  It is easy to see that this semimetric is compatible
with the usual topology on $A^*({\bf Q}_p)$, because $A^*({\bf Q}_p)$
is an open subgroup of $A^*({\bf Q}_p)$.  Let us check that
(\ref{rho(f({bf Z}_p), g({bf Z}_p))}) is also proper on $A^*({\bf Q}_p)$.

        A set $\mathcal{E} \subseteq \mathcal{C}({\bf Q}_p)$ is bounded 
with respect to $\rho(C, C')$ if and only if the cells in
$\mathcal{E}$ are contained in a single cell in ${\bf Q}_p$, and the
diameters of the cells in $\mathcal{E}$ with respect to the $p$-adic
metric on ${\bf Q}_p$ have a positive lower bound.  Of course, the
first condition implies that the diameters of the cells in
$\mathcal{E}$ also have a finite upper bound.  A set $E \subseteq
A^*({\bf Q}_p)$ is bounded with respect to (\ref{rho(f({bf Z}_p),
  g({bf Z}_p))}) if anf only if the set $\mathcal{E}$ of cells of the
form $f({\bf Z}_p)$ with $f \in E$ is bounded in $\mathcal{C}({\bf
  Q}_p)$.  If $E \subseteq A^*({\bf Q}_p)$ has this property and $f
\in E$ is as in (\ref{f(x) = a x + b}), then the coefficients $a \in
{\bf Q}_p \setminus \{0\}$ and $b \in {\bf Q}_p$ corresponding to $f$ 
are bounded in ${\bf Q}_p$, and there is a positive lower bound for
$|a|_p$.  If $E$ is also a closed set with respect to the usual
topology on $A^*({\bf Q}_p)$, then it follows that $E$ is compact, as
desired.

\end{document}